\documentclass{amsart} 
\usepackage{amsfonts}

\newtheorem{theorem}{Theorem}

\newtheorem{corollary}[theorem]{Corollary}

\newtheorem{lemma}[theorem]{Lemma}

\numberwithin{equation}{section}

\begin{document}
\title{Area and spectrum estimates for stable minimal surfaces}
\author{Ovidiu Munteanu}
\email{ovidiu.munteanu@uconn.edu}
\address{Department of Mathematics, University of Connecticut, Storrs, CT
06268, USA}
\author{Chiung-Jue Anna Sung}
\email{cjsung@math.nthu.edu.tw}
\address{Department of Mathematics, National Tsing Hua University, Hsin-Chu,
Taiwan}
\author{Jiaping Wang}
\email{jiaping@math.umn.edu}
\address{School of Mathematics, University of Minnesota, Minneapolis, MN
55455, USA}
\thanks{The first author was partially supported by NSF grant DMS-1811845
and the second by MOST, Taiwan.}

\begin{abstract}
This paper mainly concerns the area growth and bottom spectrum of complete stable
minimal surfaces in a three-dimensional manifold with scalar curvature
bounded from below. When the ambient manifold is the Euclidean space, by an
elementary argument, it is shown directly from the stability inequality that
the area of such minimal surfaces grows exactly as the Euclidean plane.
Consequently, such minimal surfaces must be flat, a well-known result due to
Fisher-Colbrie and Schoen as well as do Carmo and Peng. In the case of 
general ambient manifold, explicit area growth estimate is also derived. 
For the bottom spectrum, a self-contained argument involving positive Green's function 
is provided for its upper bound estimates. The argument extends to stable minimal hypersurfaces 
in a complete manifold of dimension up to six with sectional curvature bounded from below.
\end{abstract}

\maketitle

\section{Introduction}

This paper mainly concerns stable minimal surfaces in three-dimensional
manifolds. Our goal is to derive geometric information of such surfaces
under scalar curvature assumption on the ambient manifold. Recall that a
minimal hypersurface $\Sigma$ of a complete manifold $M$ is said to be
stable if it minimizes area up to the second order with respect to compactly
supported variations. This is equivalent to the validity of the following
stability inequality.

\begin{equation}
\int_{\Sigma }\left( \left\vert h\right\vert ^{2}+\mathrm{Ric}\left( \nu
,\nu \right) \right) \phi ^{2}\leq \int_{\Sigma }\left\vert \nabla \phi
\right\vert ^{2}  \label{S1}
\end{equation}%
for any $\phi \in C_{0}^{\infty }\left( \Sigma \right),$ where $h$ is the
second fundamental form of $\Sigma $ and $\mathrm{Ric}\left( \nu ,\nu\right)$ 
the Ricci curvature of $M$ evaluated at the unit normal vector $\nu$ to 
$\Sigma.$ As is well-known, the stability inequality is equivalent to the
existence of a positive solution $u$ to the following equation.

\begin{equation*}
\Delta u+\left( \left\vert h\right\vert ^{2}+\mathrm{Ric}\left( \nu,\nu
\right) \right)\,u=0
\end{equation*}
on $\Sigma.$ Since the lift of $u$ to the universal cover $\widetilde{\Sigma}
$ of $\Sigma$ remains a solution to the above equation, this shows that $%
\widetilde{\Sigma}$ is a stable minimal surface in the universal cover $%
\widetilde{M}$ of $M$ as well. Also, observe that such $u$ is positive
superharmonic when $M$ has nonnegative Ricci curvature. In particular, if $%
\Sigma$ is parabolic, then function $u$ must be a constant. Consequently, $%
\Sigma$ is totally geodesic.

Historically, when the ambient three-dimensional manifold $M$ has
nonnegative scalar curvature, it was shown by Schoen and Yau in their
pioneering work \cite{SY} that a compact stable minimal surface $\Sigma$
must be of genus $0.$ Subsequently, it was proven by Fisher-Colbrie and
Schoen \cite{FS} that a simply connected complete stable minimal surface $%
\Sigma$ must be conformal to the Euclidean plane, hence parabolic. As
pointed out above, this enabled them to conclude that such $\Sigma$ is
necessarily totally geodesic if the Ricci curvature of the ambient manifold $%
M$ is nonnegative. In particular, they obtained the following theorem, which
was also proved independently by do Carmo and Peng \cite{DP}.

\begin{theorem}[Fisher-Colbrie and Schoen, do Carmo and Peng]
\label{i1} A complete stable minimal surface $\Sigma$ in $\mathbb{R}^3$ must
be flat.
\end{theorem}

Later, an alternative proof of the above theorem was produced by Pogorelov \cite{P}
and Colding-Minicozzi \cite{CM1} by establishing the following area estimate

\begin{equation*}
A(r)\leq \frac{4}{3}\,\pi\,r^2
\end{equation*}
for all $r>0.$ Here, $A(r)$ denotes the area of the geodesic ball of radius 
$r$ in $\Sigma.$ It then follows \cite{CY} that $\Sigma$ must be parabolic,
hence totally geodesic and flat. In passing, we mention a recent exciting
result by Chodosh and Li \cite{CL} that a three-dimensional complete stable
minimal hypersurface in $\mathbb{R}^4$ must be flat as well.

One of our goals here is to improve the above area estimate to the sharp
form of $A(r)\leq \pi\,r^2$ for all $r>0,$ thus bypassing the parabolicity
consideration and leading directly to Theorem \ref{i1}. Indeed, since the
sectional curvature $K_{\Sigma}$ of $\Sigma$ is nonpositive and $\Sigma$ can
be assumed to be simply connected, one sees immediately from the area bound
that $\Sigma$ must be flat. In fact, our area estimate can be localized to
stable minimal surfaces with boundary. In the following, we use $B_{p}\left(
r\right)$ to denote the geodesic ball in $\Sigma$ of radius $r$ centered at
point $p\in \Sigma.$ Its area and the length of geodesic circle $\partial
B_{p}\left( r\right)$ are denoted by $A(r)$ and $L(r),$ respectively.

\begin{theorem}
\label{i2} Let $\Sigma $ be a stable minimal surface in $\mathbb{R}^{3}.$
Then there exists a universal constant $R_{0}$ such that for any geodesic
ball $B_{p}\left( R\right)$ with no intersection with the boundary of $\Sigma,$ 

\begin{equation*}
L\left( r\right) \leq 2\pi r\left( 1+\frac{10}{\ln R}\right)
\end{equation*}%
and 
\begin{equation*}
A\left( r\right) \leq \pi\, r^2\left( 1+\frac{10}{\ln R}\right)
\end{equation*}%
for all $r\leq \sqrt{R}$ and $R\geq R_{0}.$ In particular, if $\Sigma$ is
complete, then $A\left( r\right)\leq \pi\, r^2$ for all $r>0.$ Consequently, 
$\Sigma $ is flat.
\end{theorem}

More generally, we also obtain area estimate for $\Sigma$ in terms of the
lower bound of the scalar curvature $S$ of the ambient manifold $M.$

\begin{theorem}
\label{i3} Let $B_{p}\left( R\right) $ be a geodesic ball in a stable
minimal surface $\Sigma $ in a three dimensional manifold $M.$ Assume that $%
B_{p}\left( R\right) $ does not intersect the boundary of $\Sigma.$ 

\begin{itemize}
\item If the scalar curvature $S$ of $M$ satisfies $S\geq -6,$ then

\begin{equation*}
A\left( R\right) \leq C_{1}\,e^{2\,R}
\end{equation*}
for some absolute constant $C_{1}>0.$

\item If the sectional curvature $K$ of $M$ satisfies $K\geq -1,$ then

\begin{equation*}
A\left( R\right) \leq C_{1}\,e^{\frac{4}{\sqrt{7}}\,R}
\end{equation*}
for some absolute constant $C_{1}>0.$
\end{itemize}
\end{theorem}

\begin{corollary}
\label{i4} Let $\Sigma $ be a complete stable minimal surface in a three
dimensional manifold $M.$
Then there exists an absolute constant $C_{1}>0$ such that for all $R>0,$

\begin{equation*}
A\left( R\right) \leq C_{1}\,e^{\beta\,R},
\end{equation*}
where $\beta=2$ if the scalar curvature of $M$ satisfies $S\geq -6$ and 
$\beta=\frac{4}{\sqrt{7}}$ 
if the sectional curvature of $M$ satisfies $K\geq -1.$
\end{corollary}

It is unclear to us whether this estimate is sharp. It should be mentioned
that, unlike the Euclidean case, there are infinitely many non-totally
geodesic stable minimal surfaces in the hyperbolic space $\mathbb{H}^{3}.$
Rotationally symmetric examples were constructed in \cite{Mo} and \cite{DD}.

Let us now briefly indicate the main ideas behind the proofs.
Observe that for a geodesic ball $B_p(R)$ of a general surface $N,$ if it has
no intersection with either the boundary of $N$ or the cut locus of $p,$ 
then the Gauss-Bonnet formula gives

\begin{eqnarray*}
\frac{d}{dr}L(r)&=&\frac{d}{dr}\left( \int_{\partial B_{p}\left( r\right) }ds\right)\\
&=&\int_{\partial B_{p}\left( r\right) }\left( \Delta r\right) ds  
\\
&=&\int_{\partial B_{p}\left( r\right) }k_{g}ds  \\
&=&2\pi\,\chi \left( B_{p}\left( r\right) \right) -\int_{B_{p}\left( r\right) }K_{N}dA\\
&\leq& 2\pi-\int_{B_{p}\left( r\right) }K_{N}dA,
\end{eqnarray*}
where $k_{g}$ is the geodesic curvature of $\partial B_{p}\left( r\right)$ and $K_N$ the sectional
curvature of $N.$ Note that we have used the fact that the Euler characteristic number
$\chi \left( B_{p}\left( r\right) \right) \leq 1.$ 

In particular, for a complete stable minimal surface $\Sigma$ in $\mathbb{R}^3,$
since $K_{\Sigma}\leq 0,$ by working with the universal cover of $\Sigma$ if necessary, 
the inequality 

\begin{equation}
\frac{d}{dr}L(r)\leq 2\pi-\int_{B_p(r)} K_{\Sigma} \label{F}
\end{equation}
is valid for all $r.$ Combining it with the stability inequality (\ref{S1}) one concludes that 

\begin{equation}
3\int_{\Sigma}\left( \phi ^{\prime }\right)
^{2}+4\int_{\Sigma}\phi \phi ^{\prime \prime }\leq 4\pi \phi 
^{2}\left( 0\right) \label{E}
\end{equation}
holds for any compactly supported $C^{2}$ nonincreasing function $\phi =\phi \left(r\right)$ 
on $[0,\infty).$ 

So far, we have followed the argument in \cite{CM1, CM} closely. In fact, similar argument
had been adopted earlier in \cite{G}, where it was shown that the area must be of quadratic growth
for a finite index minimal surface in a three dimensional complete manifold with real analytic metric 
and nonnegative scalar curvature.

By choosing a linear cut-off function $\phi$ in (\ref{E}), one concludes immediately that 
$A(r)\leq \frac{4}{3}\,\pi\,r^2.$ For the desired conclusion that $A(r)\leq \pi\,r^2,$
a different choice of $\phi$ is needed. 

In a similar fashion, for the proof of Theorem \ref{i3}, we only use
the inequality (\ref{F}) and the stability inequality (\ref{S1}). Of course,
the preceding derivation shows that (\ref{F}) holds
for balls with no boundary points of $\Sigma$ or the cut locus of $p.$ 
It turns out that it remains true for any ball $B_p(r)$ in an arbitrary complete surface,
possibly containing cut locus of $p.$ This highly nontrivial result
is due to Fiala \cite{F} when the surface 
is the Euclidean plane endowed with a real analytic metric. Fiala's result
was later extended by Hartman \cite{H} to smooth metrics.  In its full generality, the result
is established by Shiohama and Tanaka \cite{ST1, ST2}.

Our second goal concerns upper bound of the bottom spectrum
of the Laplacian on stable minimal hypersurfaces.
The bottom spectrum of the Laplacian on a complete
manifold $N,$ denoted by $\lambda_0(N),$ is an important geometric invariant
and characterized as the optimal Poincar\'e inequality constant

\begin{equation}
\lambda _{0}\left( N \right) \int_{N}\phi ^{2}\leq \int_{N}\left\vert \nabla
\phi \right\vert ^{2} \label{P}
\end{equation}
for all compactly supported smooth function $\phi.$

According to \cite{LW}, for any $p\in N,$

\begin{equation}
\lambda _{0}\left( N \right)\leq \frac{1}{4}\,\left(\liminf_{R\to \infty}\,%
\frac{\ln V_p(R)}{R} \right)^2,  \label{V}
\end{equation}
where $V_p(R)$ denotes the volume of the geodesic ball $B_p(R)$ centered at
point $p$ of radius $R.$ 
For a stable minimal surface $\Sigma$ in three-dimensional
complete manifolds $M,$ applying Corollary \ref{i4} on the area estimate,
one immediately obtains from (\ref{V}) an upper bound estimate for $\lambda
_{0}\left( \Sigma \right),$ a result previously proven by B\'erard, Castillon, and 
Cavalcante (see Theorem 5.1 in \cite{BCC}) by a different approach.

\begin{theorem}[B\'erard, Castillon, and Cavalcante]
\label{i5} Let $\Sigma $ be a complete stable minimal surface in a three
dimensional manifold $M.$ 

\begin{itemize}
\item If the scalar curvature of $M$ satisfies $S\geq -6,$ then 
\begin{equation*}
\lambda _{0}\left( \Sigma \right) \leq 1.
\end{equation*}

\item If the sectional curvature of $M$ satisfies $K\geq -1,$ then 
\begin{equation*}
\lambda _{0}\left( \Sigma \right) \leq \frac{4}{7}.
\end{equation*}
\end{itemize}
\end{theorem}

Previously, in \cite{C}, it was shown that for complete stable minimal
surface $\Sigma$ in $\mathbb{H}^3,$ its bottom spectrum is at most $\frac{4}{3}.$ 
It would be interesting to see if the improved upper bound of $\frac{4}{7}$ is sharp.

The proof in \cite{BCC} and the one indicated above through the area estimates both rely on Fiala's inequality (\ref{F}).
To get around this rather difficult and dimension specific result, we provide yet another approach to the theorem. 
The argument instead follows the idea in \cite{LW1, MSW, MW}
and involves the minimal positive Green's function $G$. Indeed, we take $\phi=\psi\,\left\vert \nabla G\right\vert^{1/2}$
as a test function in the Poincar\'e inequality (\ref{P}), where $\psi$ is a suitably chosen cut-off function. The Bochner formula 
for $\left\vert \nabla G\right\vert$ is then used to estimate the relevant terms, with the Ricci curvature term 
arising from the formula controlled by the stability inequality (\ref{S1}). 
This approach seems to be more direct and self-contained.
Moreover, it generalizes to stable minimal hypersurfaces of dimension up to five.

\begin{theorem} \label{i6}
Let $\Sigma $ be a complete stable minimal hypersurface in $(n+1)$-dimensional
manifold $M$with $n\leq 5.$  If the sectional curvature of $M$ satisfies $K\geq -\kappa$ for some 
nonnegative constant $\kappa,$ then 

\begin{equation*}
\lambda _{0}\left( \Sigma \right) \leq \frac{2n(n-1)^2}{6n-n^2-1}\,\kappa.
\end{equation*}
\end{theorem}

Presently, it is unclear to us how to derive an upper bound for $\lambda _{0}\left( \Sigma \right)$
when $n\geq 6.$ 
It is also worth mentioning that for any complete $m$-dimensional 
minimal submanifold
$\Sigma$ in $\mathbb{H}^n,$ according to \cite{ChL}, its
bottom spectrum satisfies

\begin{equation*}
\lambda _{0}\left( \Sigma \right)\geq \frac{(m-1)^2}{4}.
\end{equation*}

The paper is arranged as follows. Section 2 is devoted to the proofs of
Theorem \ref{i2} and Theorem \ref{i3}. The estimates for the bottom spectrum
Theorem \ref{i5} and Theorem \ref{i6} are proved in Section 3.

We thank Marcos P. Cavalcante for his interest and for bringing the paper \cite{BCC} to
our attention. 
We would like to dedicate this work to Professor Peter Li on the occasion of
his seventieth birthday. All of us have benefited enormously from his
teaching and support over the years.

\section{Area estimates}

In this section, we prove both Theorem \ref{i2} and Theorem \ref{i3}. We
continue to assume that $\Sigma$ is a stable minimal surface in
three-dimensional manifold $M.$ Recall the stability inequality.

\begin{equation}
\int_{\Sigma }\left( \left\vert h\right\vert ^{2}+\mathrm{Ric}\left( \nu
,\nu \right) \right) \phi ^{2}\leq \int_{\Sigma }\left\vert \nabla \phi
\right\vert ^{2}  \label{S}
\end{equation}%
for any $\phi \in C_{0}^{\infty }\left( \Sigma \right),$ where $h$ is the
second fundamental form of $\Sigma $ and $\mathrm{Ric}\left( \nu ,\nu\right) 
$ the Ricci curvature of $M$ in the direction of the unit normal $\nu $ to $%
\Sigma.$

Fix $p\in \Sigma.$ Let 
\begin{equation*}
r\left( x\right) =d_{\Sigma }\left( p,x\right)
\end{equation*}%
be the intrinsic distance on $\Sigma $ and

\begin{equation*}
B_{p}\left( R\right) =\left\{ x\in \Sigma :r\left( x\right) <R\right\}
\end{equation*}%
the intrinsic geodesic ball of radius $R$ in $\Sigma.$ Denote with

\begin{eqnarray*}
L\left( r\right) &=&\int_{\partial B_{p}\left( r\right) }ds \\
A\left( r\right) &=&\int_{B_{p}\left( r\right) }dA
\end{eqnarray*}
the length of the geodesic circle $\partial B_{p}\left( r\right)$ and the
area of $B_{p}\left( r\right),$ respectively.

Everywhere $S,$ $\mathrm{Ric}$ and $K$ denote the scalar, Ricci and
sectional curvatures of $M,$ respectively, while $K_{\Sigma }$ denotes the
Gauss curvature of $\Sigma.$ For minimal surface $\Sigma $ in $M,$ the Gauss
curvature equation gives

\begin{equation}
K_{\Sigma }=R_{1212}-\frac{1}{2}\left\vert h\right\vert ^{2},  \label{a1}
\end{equation}%
where $\left\{e_{1}, e_{2}\right\}$ is a local orthonormal frame on $\Sigma $
and $R_{1212}$ the sectional curvature of $M$ for the two-plane spanned by $%
\left\{e_{1}, e_{2}\right\}.$ Note that since $M$ is three dimensional,

\begin{equation*}
R_{1212}+\mathrm{Ric}\left( \nu ,\nu \right) =\frac{1}{2}S.
\end{equation*}%
Therefore,

\begin{equation}
K_{\Sigma }=\frac{1}{2}S-\left( \mathrm{Ric}\left( \nu ,\nu \right) +\frac{1%
}{2}\left\vert h\right\vert ^{2}\right).  \label{a2}
\end{equation}

In the following, we present detailed proofs for Theorem \ref{i2} and Theorem \ref{i3} by assuming
that the ball $B_p(R)$ contains no cut locus of $p$ in $\Sigma.$ This suffices for proving Theorem \ref{i2}
in its full generality. It also proves Theorem \ref{i3} in the case that $M$ has nonpositive sectional
curvature. Indeed, by (\ref{a1}), the sectional curvature $K_{\Sigma }\leq 0 $ when $R_{1212}\leq 0.$ 
By working on the universal covering space of $\Sigma$ if necessary, geodesic ball
$B_p(R)$ is free of cut locus for all $R>0$ as $\Sigma$ is a Cartan-Hadamard manifold. 

The same argument also applies to the general case of Theorem \ref{i3} by invoking Fiala's inequality (\ref{F}) 
alluded in the introduction.

We start with the following lemma which has appeared in \cite{G, CM1} or Theorem 2.8 in \cite{CM},
as well as \cite{M, Ca}.

\begin{lemma}
\label{M} Let $\Sigma $ be stable minimal surface in a three dimensional
manifold $M.$ Let $B_{p}\left( R\right) $ be a geodesic ball in $\Sigma $
that does not intersect the cut locus of $p$ in $\Sigma $ or the boundary of 
$\Sigma.$ Assume that $\phi=\phi \left( r\right) $ is a Lipschitz
continuous, nonincreasing function on $\left[ 0,R\right] $ with $\phi \left(
R\right) =0.$

\begin{itemize}
\item If the scalar curvature of $M$ satisfies $S\geq -6\alpha$ for some $%
\alpha\geq 0,$ then

\begin{equation}
-2\int_{0}^{R}\phi \left( r\right) \phi ^{\prime }\left( r\right) L^{\prime
}\left( r\right) dr\leq 2\pi \phi ^{2}\left( 0\right) +\int_{B_{p}\left(
R\right) }\left( \phi ^{\prime }\right) ^{2}+3\alpha\,\int_{B_{p}\left(
R\right) }\phi ^{2}.  \label{M1}
\end{equation}

\item If the sectional curvature of $M$ satisfies $K\geq -\alpha$ for some $%
\alpha\geq 0,$ then

\begin{equation}
-4\int_{0}^{R}\phi \left( r\right) \phi ^{\prime }\left( r\right) L^{\prime
}\left( r\right) dr\leq 4\pi \phi ^{2}\left( 0\right) +\int_{B_{p}\left(
R\right) }\left( \phi ^{\prime }\right) ^{2}+4\alpha\,\int_{B_{p}\left(
R\right) }\phi ^{2}.  \label{M2}
\end{equation}
\end{itemize}
\end{lemma}

\begin{proof}
For any $0<r<R,$ by the Gauss-Bonnet formula we have

\begin{eqnarray}
\frac{d}{dr}L(r)&=&\frac{d}{dr}\left( \int_{\partial B_{p}\left( r\right) }ds\right)  \label{GB}\\
&=&\int_{\partial B_{p}\left( r\right) }\left( \Delta r\right) ds \notag \\
&=&\int_{\partial B_{p}\left( r\right) }k_{g}ds  \notag \\
&=&2\pi \chi \left( B_{p}\left( r\right) \right)-\int_{B_{p}\left( r\right) }K_{\Sigma }dA  \notag \\
&\leq & 2\pi-\int_{B_{p}\left( r\right) }K_{\Sigma }dA,  \notag
\end{eqnarray}%
where $k_{g}$ is the geodesic curvature of $\partial B_{p}\left( r\right).$
Note that the Euler characteristic number $\chi \left( B_{p}\left( r\right) \right) \leq 1.$

(a) Assume first that $S\geq -6\alpha$. By (\ref{a2}) we have that 
\begin{equation*}
K_{\Sigma }\geq -3\alpha-\left( \mathrm{Ric}\left( \nu ,\nu \right)
+\left\vert h\right\vert ^{2}\right).
\end{equation*}
Hence,

\begin{equation*}
\frac{d}{dr}\left( \int_{\partial B_{p}\left( r\right) }ds\right) \leq 2\pi
+3\alpha\,\int_{B_{p}\left( r\right) }dA+\int_{B_{p}\left( r\right) }\left( 
\mathrm{Ric}\left( \nu ,\nu \right) +\left\vert h\right\vert ^{2}\right).
\end{equation*}%
Equivalently,

\begin{equation}
L^{\prime }\left( r\right) \leq 2\pi +3\alpha\,A\left( r\right)
+\int_{B_{p}\left( r\right) }\left( \mathrm{Ric}\left( \nu ,\nu \right)
+\left\vert h\right\vert ^{2}\right).  \label{a3}
\end{equation}%
Multiply (\ref{a3}) by $-2\phi \left( r\right) \phi ^{\prime }\left(r\right)
\geq 0$ and integrate from $r=0$ to $r=R.$ It follows that

\begin{eqnarray}
&&-2\int_{0}^{R}\phi \left( r\right) \phi ^{\prime }\left( r\right)
L^{\prime }\left( r\right) dr  \label{a4} \\
&\leq &2\pi \phi ^{2}\left( 0\right) -6\alpha\,\int_{0}^{R}\phi \left(
r\right) \phi ^{\prime }\left( r\right) A\left( r\right) dr  \notag \\
&&-2\int_{0}^{R}\phi \left( r\right) \phi ^{\prime }\left( r\right) \left(
\int_{B_{p}\left( r\right) }\mathrm{Ric}\left( \nu ,\nu \right) +\left\vert
h\right\vert ^{2}\right) dr.  \notag
\end{eqnarray}%
Note that for any function $f\left( r\right) $ with $f\left( 0\right) =0$ we
have 
\begin{equation}
-2\int_{0}^{R}\phi \left( r\right) \phi ^{\prime }\left( r\right) f\left(
r\right) dr=\int_{0}^{R}f^{\prime }\left( r\right) \phi ^{2}\left( r\right)
dr.  \label{a5}
\end{equation}%
Applying this to $f\left( r\right) =\int_{B_{p}\left( r\right) }dA$ we
conclude

\begin{eqnarray}
-6\alpha\,\int_{0}^{R}\phi \left( r\right) \phi ^{\prime }\left( r\right)
A\left( r\right) dr &=&3\alpha\,\int_{0}^{R}\phi ^{2}\left( r\right) L\left(
r\right) dr  \label{a5a} \\
&=&3\alpha\,\int_{B_{p}\left( R\right) }\phi ^{2},  \notag
\end{eqnarray}%
where in the last line we have used the co-area formula. Similarly, for

\begin{equation*}
f\left(r\right) =\int_{B_{p}\left( r\right) }\left( \mathrm{Ric}\left( \nu
,\nu \right) +\left\vert h\right\vert ^{2}\right)
\end{equation*}
we get

\begin{eqnarray}
&&-2\int_{0}^{R}\phi \left( r\right) \phi ^{\prime }\left( r\right) \left(
\int_{B_{p}\left( r\right) }\mathrm{Ric}\left( \nu ,\nu \right) +\left\vert
h\right\vert ^{2}\right)  \label{a5b} \\
&=&\int_{0}^{R}\phi ^{2}\left( r\right) \left( \int_{\partial B_{p}\left(
r\right) }\mathrm{Ric}\left( \nu ,\nu \right) +\left\vert h\right\vert
^{2}\right) dr  \notag \\
&=&\int_{B_{p}\left( R\right) }\left( \mathrm{Ric}\left( \nu ,\nu \right)
+\left\vert h\right\vert ^{2}\right) \phi ^{2}  \notag \\
&\leq &\int_{B_{p}\left( R\right) }\left( \phi ^{\prime }\right) ^{2}, 
\notag
\end{eqnarray}%
where in the last line we have used the stability inequality (\ref{S}).
Combining (\ref{a5a}) and (\ref{a5b}), we conclude from (\ref{a4}) that 
\begin{equation*}
-2\int_{0}^{R}\phi \left( r\right) \phi ^{\prime }\left( r\right) L^{\prime
}\left( r\right) dr\leq 2\pi \phi ^{2}\left( 0\right)
+3\alpha\,\int_{B_{p}\left( R\right) }\phi ^{2}+\int_{B_{p}\left( R\right)
}\left( \phi ^{\prime }\right) ^{2},
\end{equation*}
which is (\ref{M1}).

(b) Assume now that $K\geq -\alpha.$ Then according to (\ref{a1}),

\begin{equation*}
-2K_{\Sigma }\leq 2\alpha+\left\vert h\right\vert ^{2}.
\end{equation*}%
Plugging into (\ref{GB}) gives

\begin{equation*}
2L^{\prime }\left( r\right) \leq 4\pi +2\alpha\,A\left( r\right)
+\int_{B_{p}\left( r\right) }\left\vert h\right\vert ^{2}.
\end{equation*}%
Multiplying by $-2\phi \left( r\right) \phi ^{\prime }\left( r\right)\geq 0$
and integrating from $r=0$ to $r=R$ we obtain

\begin{eqnarray}
-4\int_{0}^{R}\phi \left( r\right) \phi ^{\prime }\left( r\right) L^{\prime
}\left( r\right) dr &\leq &4\pi \phi ^{2}\left( 0\right)  \label{a7} \\
&&-4\alpha\,\int_{0}^{R}\phi \left( r\right) \phi ^{\prime }\left( r\right)
A\left( r\right) dr  \notag \\
&&-2\int_{0}^{R}\phi \left( r\right) \phi ^{\prime }\left( r\right) \left(
\int_{B_{p}\left( r\right) }\left\vert h\right\vert ^{2}\right) dr.  \notag
\end{eqnarray}
Since $\mathrm{Ric}\left( \nu ,\nu \right)\geq -2\alpha,$ it follows from (%
\ref{a5b}) that

\begin{equation*}
-2\int_{0}^{R}\phi \left( r\right) \phi ^{\prime }\left( r\right) \left(
\int_{B_{p}\left( r\right) }\left\vert h\right\vert ^{2}\right) \leq
\int_{\Sigma }\left( \phi ^{\prime }\right) ^{2}+2\alpha\,\int_{\Sigma }\phi
^{2}.
\end{equation*}
Together with (\ref{a5a}), we conclude from (\ref{a7}) that

\begin{equation*}
-4\int_{0}^{R}\phi \left( r\right) \phi ^{\prime }\left( r\right) L^{\prime
}\left( r\right) dr\leq 4\pi \,\phi ^{2}\left( 0\right) +\int_{B_{p}\left(
R\right) }\left( \phi ^{\prime }\right) ^{2}+4\alpha\,\int_{B_{p}\left(
R\right) }\phi ^{2}.
\end{equation*}
This proves (\ref{M2}).
\end{proof}

\begin{corollary}
\label{C1} Let $\Sigma $ be stable minimal surface in a three dimensional
manifold $M.$ Let $B_{p}\left( R\right) $ be a geodesic ball in $\Sigma $
that does not intersect the cut locus of $p$ in $\Sigma $ or the boundary of 
$\Sigma.$ Assume that $\phi=\phi \left( r\right) $ is a $C^{2}$
nonincreasing function on $\left[ 0,R\right] $ with $\phi \left( R\right)
=0. $

\begin{itemize}
\item If the scalar curvature of $M$ satisfies $S\geq -6\alpha$ for some $%
\alpha\geq 0,$ then 
\begin{equation*}
\int_{B_{p}\left( R\right) }\left( \phi ^{\prime }\right)
^{2}+2\int_{B_{p}\left( R\right) }\phi \phi ^{\prime \prime }\leq 2\pi \phi
^{2}\left( 0\right) +3\alpha\,\int_{B_{p}\left( R\right) }\phi ^{2}.
\end{equation*}

\item If the sectional curvature of $M$ satisfies $K\geq -\alpha$ for some $%
\alpha\geq 0,$ then 
\begin{equation*}
3\int_{B_{p}\left( R\right) }\left( \phi ^{\prime }\right)
^{2}+4\int_{B_{p}\left( R\right) }\phi \phi ^{\prime \prime }\leq 4\pi \phi
^{2}\left( 0\right) +4\alpha\,\int_{B_{p}\left( R\right) }\phi ^{2}.
\end{equation*}
\end{itemize}
\end{corollary}

\begin{proof}
Integrating by parts we have

\begin{eqnarray*}
-2\int_{0}^{R}\phi \left( r\right) \phi ^{\prime }\left( r\right) L^{\prime
}\left( r\right) dr &=&2\int_{0}^{R}\left( \phi \phi ^{\prime \prime
}+\left( \phi ^{\prime }\right) ^{2}\right) L\left( r\right) dr \\
&=&2\int_{B_{p}\left( R\right) }\left( \phi \phi ^{\prime \prime }+\left(
\phi ^{\prime }\right) ^{2}\right).
\end{eqnarray*}%
By (\ref{M1}) and (\ref{M2}), the desired conclusions follow.
\end{proof}

We now use the lemma and the corollary to establish area estimates for
stable minimal surfaces. We start with the case $M=\mathbb{R}^{3}.$

\begin{theorem}
\label{v2} Let $\Sigma $ be a stable minimal surface in $\mathbb{R}^{3}.$
Then there exists a universal constant $R_{0}$ such that for geodesic ball $%
B_{p}\left( R\right)$ with no intersection with the boundary of $\Sigma$ or
the cut locus of $p,$

\begin{equation*}
L\left( r\right) \leq 2\pi r\left( 1+\frac{10}{\ln R}\right)
\end{equation*}%
and 
\begin{equation*}
A\left( r\right) \leq \pi\, r^2\left( 1+\frac{10}{\ln R}\right)
\end{equation*}%
for all $r\leq \sqrt{R}$ and $R\geq R_{0}.$ In particular, if $\Sigma$ is
complete, then $A\left( r\right)\leq \pi\, r^2$ for all $r>0.$ Consequently, 
$\Sigma $ is flat.
\end{theorem}

\begin{proof}
By Corollary \ref{C1}, the inequality

\begin{equation*}
3\int_{B_{p}\left( R\right) }\left( \phi ^{\prime }\right)
^{2}+4\int_{B_{p}\left( R\right) }\phi \phi ^{\prime \prime }\leq 4\pi \phi
^{2}\left( 0\right)
\end{equation*}
holds for any $C^{2}$ nonincreasing function $\phi =\phi \left(r\right)$ on $%
\left[ 0,R\right] $ with $\phi \left(R\right) =0.$ Set

\begin{equation*}
\phi \left( r\right) =\ln \left( R+1\right) -\ln \left( r+1\right).
\end{equation*}
The inequality becomes

\begin{equation*}
3\int_{B_{p}\left( R\right) }\frac{1}{\left( r+1\right) ^{2}}
+4\int_{B_{p}\left( R\right) }\frac{\ln \left( R+1\right) -\ln \left(
r+1\right) }{\left( r+1\right) ^{2}}\leq 4\pi \ln ^{2}\left( R+1\right).
\end{equation*}
In particular,

\begin{equation}
\int_{B_{p}\left( R\right) }\frac{\ln \left( R+1\right) -\ln \left(
r+1\right) }{\left( r+1\right) ^{2}}\leq \pi \ln ^{2}\left( R+1\right).
\label{a16}
\end{equation}

Since $K_{\Sigma }\leq 0$ on $B_{p}\left( R\right)$, the Hessian comparison
theorem (see Theorem 1.1 in \cite{SY1}) implies that

\begin{equation}
2\pi \leq \frac{L\left( r\right) }{r}\leq \frac{L\left( R\right) }{R}
\label{a17}
\end{equation}%
for all $0<r<R.$ Assume by contradiction that

\begin{equation}
\frac{L\left( r\right) }{r}\geq 2\pi \left( 1+\frac{10}{\ln \left(
R+1\right) }\right)  \label{a18}
\end{equation}
for all $r\in \left[ \sqrt{R},R\right].$ According to (\ref{a16}) we have

\begin{eqnarray*}
\pi \ln ^{2}\left( R+1\right) &\geq &\int_{0}^{R}\frac{\ln \left( R+1\right)
-\ln \left( r+1\right) }{\left( r+1\right) ^{2}}L\left( r\right) dr \\
&=&\int_{0}^{R}\frac{\ln \left( R+1\right) -\ln \left( r+1\right) }{\left(
r+1\right) ^{2}}\left( L\left( r\right) -2\pi r\right) dr \\
&&+2\pi \int_{0}^{R}\frac{\ln \left( R+1\right) -\ln \left( r+1\right) }{%
\left( r+1\right) ^{2}}r\,dr.
\end{eqnarray*}
Hence, (\ref{a17}) and (\ref{a18}) imply that 
\begin{eqnarray}
\pi \ln ^{2}\left( R+1\right) &\geq &2\pi \int_{0}^{R}\frac{\ln \left(
R+1\right) -\ln \left( r+1\right) }{\left( r+1\right) ^{2}}r\,dr  \label{a19}
\\
&&+\frac{20\pi }{\ln \left( R+1\right) }\int_{\sqrt{R}}^{R}\frac{\ln \left(
R+1\right) -\ln \left( r+1\right) }{\left( r+1\right) ^{2}}rdr.  \notag
\end{eqnarray}%
The first term can be computed as

\begin{eqnarray*}
&&2\pi \int_{0}^{R}\frac{\ln \left( R+1\right) -\ln \left( r+1\right) }{%
\left( r+1\right) ^{2}}rdr \\
&=&-2\pi \ln \left( R+1\right) +2\pi \int_{0}^{R}\frac{1}{r+1}\left( \ln
\left( r+1\right) +\frac{1}{r+1}\right) dr \\
&=&-2\pi \ln \left( R+1\right) +\pi \ln ^{2}\left( R+1\right) +2\pi -\frac{%
2\pi }{R+1}.
\end{eqnarray*}%
Plugging this into (\ref{a19}) we get

\begin{equation}
2\pi \ln \left( R+1\right) \geq \frac{20\pi }{\ln \left( R+1\right) }\int_{%
\sqrt{R}}^{R}\frac{\ln \left( R+1\right) -\ln \left( r+1\right) }{\left(
r+1\right) ^{2}}r\,dr.  \label{a20}
\end{equation}%
However, integration by parts gives

\begin{eqnarray*}
&&\int_{\sqrt{R}}^{R}\frac{\ln \left( R+1\right) -\ln \left( r+1\right) }{%
\left( r+1\right) ^{2}}rdr \\
&=&-\ln \left( \frac{R+1}{\sqrt{R}+1}\right) \left( \ln \left( \sqrt{R}%
+1\right) +\frac{1}{\sqrt{R}+1}\right) \\
&&+\frac{1}{2}\ln ^{2}\left( R+1\right) -\frac{1}{2}\ln ^{2}\left( \sqrt{R}%
+1\right) \\
&&+\frac{1}{\sqrt{R}+1}-\frac{1}{R+1} \\
&\geq &\frac{1}{9}\ln ^{2}\left( R+1\right)
\end{eqnarray*}%
for all $R>R_{0}$ large enough. In view of (\ref{a20}), this yields a
contradiction. In conclusion, (\ref{a18}) is false. In other words, there
exists $r_{0}\in \left[ \sqrt{R},R\right] $ such that

\begin{equation*}
\frac{L\left( r_{0}\right) }{r_{0}}\leq 2\pi \left( 1+\frac{10}{\ln \left(
R+1\right) }\right).
\end{equation*}
Now for $r<\sqrt{R},$ by (\ref{a17}) we have

\begin{eqnarray*}
\frac{L\left( r\right) }{r} &\leq &\frac{L\left( r_{0}\right) }{r_{0}} \\
&\leq &2\pi \left( 1+\frac{10}{\ln R}\right).
\end{eqnarray*}%
This proves the length estimate. The area estimate then follows immediately.

Finally, if $\Sigma$ is complete, then the length and area estimates are
applicable for all $R$ on the universal cover of $\Sigma.$ It follows that

\begin{equation*}
L\left( r\right)\leq 2\pi\,r
\end{equation*}
and 
\begin{equation*}
A\left( r\right)\leq \pi\,r^2
\end{equation*}
for all $r>0.$ This implies that $\Sigma $ is flat as $K_{\Sigma}\leq 0.$
\end{proof}

We now turn to the case of more general ambient manifolds.

\begin{theorem}
\label{v1} Let $B_{p}\left( R\right) $ be a geodesic ball in a stable
minimal surface $\Sigma $ in a three dimensional manifold $M.$ Assume that $%
B_{p}\left( R\right) $ does not intersect the boundary of $\Sigma$ or the
cut locus of $p$ in $\Sigma.$

\begin{itemize}
\item If the scalar curvature $S$ of $M$ satisfies $S\geq -6,$ then

\begin{equation*}
A\left( R\right) \leq C_{1}\,e^{2\,R}
\end{equation*}
for some absolute constant $C_{1}>0.$

\item If the sectional curvature $K$ of $M$ satisfies $K\geq -1,$ then

\begin{equation*}
A\left( R\right) \leq C_{1}\,e^{\frac{4}{\sqrt{7}}\,R}
\end{equation*}
for some absolute constant $C_{1}>0.$
\end{itemize}
\end{theorem}

\begin{proof}
Since the arguments are similar for both cases, we supply details for the
second case and only provide a sketch for the first case.

So we assume $K\geq -1.$ By (\ref{M2}) from Lemma \ref{M}, for any Lipschitz
continuous, nonincreasing function $\phi =\phi \left( r\right) $ on $\left[
0,t\right] $ with $\phi \left(t\right) =0$ we have

\begin{equation}
-4\int_{0}^{t}\phi \left( r\right) \phi ^{\prime }\left( r\right) L^{\prime
}\left( r\right) dr\leq 4\pi \phi ^{2}\left( 0\right) +\int_{B_{p}\left(
t\right) }\left( \phi ^{\prime }\right) ^{2}+4\int_{B_{p}\left( t\right)
}\phi ^{2}  \label{a8}
\end{equation}%
for all $0<t<R.$

For convenience, denote with 
\begin{equation}
a=\frac{4}{\sqrt{7}}.  \label{a}
\end{equation}
Let 
\begin{equation*}
\phi \left( r\right) =e^{-\frac{a}{2}r}\psi \left( r\right),
\end{equation*}%
where $\psi $ is a nonincreasing Lipschitz function such that $\psi
\left(t\right) =0.$ We have

\begin{eqnarray*}
\left( \phi ^{\prime }\right) ^{2} &=&\frac{a^{2}}{4}e^{-ar}\psi
^{2}+e^{-ar}\left( \psi ^{\prime }\right) ^{2}-ae^{-ar}\psi \psi ^{\prime }
\\
-4\phi \phi ^{\prime } &=&2ae^{-ar}\psi ^{2}-4e^{-ar}\psi \psi ^{\prime }.
\end{eqnarray*}
Therefore,

\begin{eqnarray*}
-4\int_{0}^{t}\phi \left( r\right) \phi ^{\prime }\left( r\right) L^{\prime
}\left( r\right) dr &=&2a\int_{0}^{t}e^{-ar}\psi ^{2}\left( r\right)
L^{\prime }\left( r\right) dr \\
&&-4\int_{0}^{t}e^{-ar}\psi \left( r\right) \psi ^{\prime }\left( r\right)
L^{\prime }\left( r\right) dr.
\end{eqnarray*}%
After integration by parts, the first term on the right side becomes

\begin{eqnarray*}
2a\int_{0}^{t}e^{-ar}\psi ^{2}\left( r\right) L^{\prime }\left( r\right) dr
&=&\int_{0}^{t}\left( 2a^{2}e^{-ar}\psi ^{2}\left( r\right) -4a\psi \left(
r\right) \psi ^{\prime }\left( r\right) \right) L\left( r\right) dr \\
&=&\int_{B_{p}\left( t\right) }\left( 2a^{2}\psi ^{2}-4a\psi \psi ^{\prime
}\right) e^{-ar}.
\end{eqnarray*}

Plugging these identities into (\ref{a8}), we get

\begin{eqnarray}
-4\int_{0}^{t}e^{-ar}\psi \left( r\right) \psi ^{\prime }\left( r\right)
L^{\prime }\left( r\right) dr &\leq &4\pi \psi ^{2}\left( 0\right)
+\int_{B_{p}\left( t\right) }\left( \psi ^{\prime }\right) ^{2}e^{-ar}
\label{a9} \\
&&+3a\int_{B_{p}\left( t\right) }\psi \psi ^{\prime }e^{-ar}.  \notag
\end{eqnarray}

For arbitrary $\eta$ with $0<\eta<R$ and $\eta\leq t$, let

\begin{equation}
\psi \left( r\right) =\left\{ 
\begin{array}{c}
1 \\ 
\frac{t-r}{\eta } \\ 
0%
\end{array}%
\right. 
\begin{array}{l}
\text{for }r\leq t-\eta \\ 
\text{for }r\in \left( t-\eta ,t\right) \\ 
\text{for }r\geq t%
\end{array}
\label{psi}
\end{equation}%
It follows from (\ref{a9}) that

\begin{eqnarray}
\frac{4}{\eta }\int_{t-\eta }^{t}e^{-ar}\psi \left( r\right) L^{\prime
}\left( r\right) dr &\leq &4\pi +\frac{1}{\eta ^{2}}\int_{B_{p}\left(
t\right) \backslash B_{p}\left( t-\eta \right) }e^{-ar}  \label{a10} \\
&&-\frac{3a}{\eta }\int_{B_{p}\left( t\right) \backslash B_{p}\left( t-\eta
\right) }\psi e^{-ar}.  \notag
\end{eqnarray}%
The term on the left hand side, after integrating by parts, becomes

\begin{eqnarray*}
\frac{4}{\eta }\int_{t-\eta }^{t}e^{-ar}\psi \left( r\right) L^{\prime
}\left( r\right) dr &=&-\frac{4}{\eta }L\left( t-\eta \right) e^{-a\left(
t-\eta \right) }+\frac{4a}{\eta }\int_{B_{p}\left( t\right) \backslash
B_{p}\left( t-\eta \right) }\psi e^{-ar} \\
&&+\frac{4}{\eta ^{2}}\int_{B_{p}\left( t\right) \backslash B_{p}\left(
t-\eta \right) }e^{-ar}.
\end{eqnarray*}%
Therefore, combining with (\ref{a10}), we have

\begin{equation}
\frac{3}{\eta ^{2}}\int_{B_{p}\left( t\right) \backslash B_{p}\left( t-\eta
\right) }e^{-ar}+\frac{7a}{\eta }\int_{B_{p}\left( t\right) \backslash
B_{p}\left( t-\eta \right) }\psi e^{-ar}\leq 4\pi +\frac{4}{\eta }L\left(
t-\eta \right) e^{-a\left( t-\eta \right) }  \label{a11}
\end{equation}%
for any $0<\eta \leq t<R.$ In particular, by setting $\eta=t,$ it follows
that

\begin{equation*}
\frac{3}{t^{2}}\int_{B_{p}\left( t\right) }e^{-ar}\leq 4\pi.
\end{equation*}%
Clearly, this implies that

\begin{equation}
A\left( t\right) \leq \frac{4\pi }{3}\,t^{2}\,e^{at}  \label{a11'}
\end{equation}%
for all $t\leq R.$

\vspace{0.1 in}

\noindent \textbf{Claim:} There exists an absolute constant $\Lambda>0$ such
that for all $\tau$ and $s$ satisfying $0<2\tau \leq s<R-3\tau$,

\begin{equation}
\int_{B_{p}\left( s\right) \backslash B_{p}\left( s-\tau \right)}e^{-ar}\leq
\Lambda\tau +\frac{\Lambda}{\tau }\int_{B_{p}\left( s-\tau \right)
\backslash B_{p}\left( s-2\tau \right) }e^{-ar}.  \label{a14}
\end{equation}

Indeed, let $\eta=4\tau$ and $T=s-\frac{3\tau}{2}.$ Then

\begin{equation*}
\frac{\eta }{8}\leq T< R-\frac{9\eta }{8}.
\end{equation*}
The mean value theorem implies that there exists $\xi \in \left( T-\frac{\eta%
}{8},T+\frac{\eta }{8}\right)$ such that

\begin{equation}
\int_{B_{p}\left( T+\frac{\eta}{8}\right) \backslash B_{p}\left( T-\frac{\eta%
}{8}\right) }e^{-ar} =\frac{\eta }{4}L\left( \xi \right) e^{-a\xi }.
\label{a12}
\end{equation}
Denote with 
\begin{equation}
t=\xi +\eta \in \left( T+\frac{7\eta }{8},T+\frac{9\eta }{8}\right).
\label{t}
\end{equation}
By (\ref{psi}) and (\ref{t}) we get that

\begin{eqnarray*}
\int_{B_{p}\left( t\right) \backslash B_{p}\left( t-\eta \right) }\psi
e^{-ar} &\geq &\frac{1}{2}\int_{B_{p}\left( t-\frac{\eta}{2}\right)
\backslash B_{p}\left( t-\eta \right) }e^{-ar} \\
&\geq &\frac{1}{2}\int_{B_{p}\left( T+\frac{3\eta }{8}\right) \backslash
B_{p}\left( T+\frac{\eta }{8}\right) }e^{-ar}.
\end{eqnarray*}%
Together with (\ref{a11}) and (\ref{a12}) we conclude

\begin{eqnarray*}
\frac{7a}{2\eta }\int_{B_{p}\left( T+\frac{3\eta }{8}\right) \backslash
B_{p}\left( T+\frac{\eta }{8}\right) }e^{-ar} &\leq &\frac{7a}{\eta }%
\int_{B_{p}\left( t\right) \backslash B_{p}\left( t-\eta \right) }\psi
e^{-ar} \\
&\leq &4\pi +\frac{4}{\eta }L\left( t-\eta \right) e^{-a\left( t-\eta
\right) } \\
&=&4\pi +\frac{16}{\eta^{2}}\int_{B_{p}\left( T+\frac{\eta }{8}\right)
\backslash B_{p}\left( T-\frac{\eta }{8}\right) }e^{-ar}.
\end{eqnarray*}%
Therefore, there exists $\Lambda>0$ such that

\begin{equation}
\int_{B_{p}\left( T+\frac{3\eta }{8}\right) \backslash B_{p}\left( T+\frac{%
\eta}{8}\right) }e^{-ar}\leq \Lambda\eta +\frac{\Lambda}{\eta }%
\int_{B_{p}\left( T+\frac{\eta }{8}\right) \backslash B_{p}\left( T-\frac{%
\eta }{8}\right) }e^{-ar}.  \label{a13}
\end{equation}%
Substituting $\eta =4\tau $ and $s=T+\frac{3\tau }{2}$ in (\ref{a13})
implies the claim.

For $s\geq 6\Lambda,$ letting $\tau=2\Lambda$ and iterating (\ref{a14}) $m$
times with $m=\left[ \frac{s}{2\Lambda}\right] -2$ we get

\begin{eqnarray}
\int_{B_{p}\left( s\right) \backslash B_{p}\left( s-2\Lambda \right)}e^{-ar}
&\leq& 2\Lambda^2\,\sum_{k=0}^{m-1}\frac{1}{2^{k}} +\frac{1}{2^{m}}%
\int_{B_{p}\left(6\Lambda\right)}e^{-ar}  \label{a15} \\
&\leq& C_{2}  \notag
\end{eqnarray}
by invoking (\ref{a11'}) for the last term. This holds for any $6\Lambda\leq
s \leq R-6\Lambda$.

To finish, we apply the mean value theorem to conclude that there exists $%
\xi \in \left( R-16\Lambda ,R-14\Lambda\right) $ such that

\begin{equation}
\int_{B_{p}\left( R-14\Lambda \right) \backslash B_{p}\left(
R-16\Lambda\right) }e^{-ar}=2\Lambda L\left( \xi \right) e^{-a\xi }.
\label{a15'}
\end{equation}
Applying (\ref{a15}) with $s=R-14\Lambda,$ we get from (\ref{a15'}) that

\begin{equation*}
L\left( R-\eta\right) e^{-a\left( R-\eta \right) }\leq C_{3}
\end{equation*}%
for some constant $C_3,$ where 
\begin{equation*}
\eta =R-\xi \in \left(14\Lambda, 16\Lambda\right).
\end{equation*}

By (\ref{a11}),

\begin{eqnarray*}
\frac{3}{\eta^{2}}\int_{B_{p}\left( R\right) \backslash B_{p}\left( R-\eta
\right) }e^{-ar} &\leq &4\pi +\frac{4}{\eta}L\left( R-\eta \right)
e^{-a\left( R-\eta \right) } \\
&\leq &C_{4}.
\end{eqnarray*}%
This implies that

\begin{equation*}
\int_{B_{p}\left( R\right) \backslash B_{p}\left( R-\eta \right) }dA\leq
C\,e^{a\,R}
\end{equation*}%
for some $\eta \in \left(14\Lambda, 16\Lambda\right).$ In particular,

\begin{equation}
\int_{B_{p}\left( R\right) \backslash B_{p}\left( R-14\Lambda\right) }dA\leq
C\,e^{a\,R}.  \label{A}
\end{equation}
Applying (\ref{A}) with $R$ replaced by $R-14k\Lambda,$ $k=1,2,\cdots, n$
and $n=\left[ \frac{R}{14\Lambda}\right] -1,$ we arrive at

\begin{eqnarray*}
A(R)&\leq& \sum_{k=0}^n \int_{B_{p}\left( R-14k\Lambda\right) \backslash
B_{p}\left( R-14(k+1)\Lambda\right) }dA+A(14\Lambda) \\
&\leq& C\,\sum_{k=0}^n e^{a\,\left(R-14k\Lambda\right)}+C \\
&\leq& C\, e^{a\,R}.
\end{eqnarray*}
This proves the area estimate for the second case.

Let us now sketch the argument for the first case. Assume now that $S\geq
-6. $ Then by (\ref{M1}),

\begin{equation*}
-2\int_{0}^{t}\phi \left( r\right) \phi ^{\prime }\left( r\right) L^{\prime
}\left( r\right) dr\leq 2\pi \phi ^{2}\left( 0\right) +\int_{B_{p}\left(
t\right) }\left( \phi ^{\prime }\right) ^{2}+3\int_{B_{p}\left( t\right)
}\phi ^{2}.
\end{equation*}
We set 
\begin{equation*}
\phi \left( r\right) =e^{-r}\psi \left( r\right),
\end{equation*}
for $\psi$ nonincreasing Lipschitz on $\left[ 0,t \right ]$ so that $\psi
(t)=0$, and get that

\begin{equation*}
-2\int_{0}^{t}\psi \left( r\right) \psi ^{\prime }\left( r\right) L^{\prime
}\left( r\right) e^{-2r}\leq 2\pi +2\int_{B_{p}\left( t\right) }\psi \psi
^{\prime }e^{-2r}+\int_{B_{p}\left( t\right) }\left( \psi ^{\prime }\right)
^{2}e^{-2r}.
\end{equation*}%
Taking $\psi $ as defined in (\ref{psi}) and integrating by parts imply that 
\begin{equation*}
\frac{1}{\eta ^{2}}\int_{B_{p}\left( t\right) \backslash B_{p}\left( t-\eta
\right) }e^{-2r}+\frac{6}{\eta }\int_{B_{p}\left( t\right) \backslash
B_{p}\left( t-\eta \right) }\psi e^{-2r}\leq 2\pi +\frac{2}{\eta}L\left(
t-\eta \right) e^{-2\left( t-\eta\right) }
\end{equation*}
The rest of the argument follows verbatim.
\end{proof}

We point out that the same proof applies to the general case of Theorem \ref{i3} as well, 
that is, Theorem \ref{v1} continues to hold for arbitrary ball $B_p(R)$ in $\Sigma.$
Indeed,  Lemma \ref{M} is valid for arbitrary ball $B_p(R),$ even if 
it contains cut locus of $p,$ by invoking Fiala's inequality (\ref{F})
established through the work of Fiala \cite{F}, Hartman \cite{H},
and Shiohama and Tanaka \cite{ST1, ST2}.
A priori, the length function $L(r)$ of the geodesic circle $\partial B_p(r)$ 
is only defined for almost all $r.$
Their work implies that the function $L(r)$ can be extended to all $r\geq 0.$  
So extended function, denoted by $L$ again, may not be continuous. But 
it satisfies $L(r^{+}) = L(r)$ and $L(r^{-})\geq L(r)$ for all $r>0.$ Moreover,
it can be written as $L(r)=L_1(r)+L_2(r),$ where
$L_1(r)$ is absolutely continuous on any finite interval and $L_2(r)$ is 
nonincreasing. At the points where $L(r)$ is differentiable, its derivative
satisfies the Fiala inequality

\begin{equation*}
\frac{d}{dr}L(r)\leq 2\pi-\int_{B_p(r)} K_{\Sigma}.
\end{equation*}

Although it is not needed here, we remark that Corollary \ref{C1} is also valid on arbitrary ball
$B_p(R)$ for all $C^{2}$ function $\phi=\phi \left( r\right) $ with $\phi'(r)\leq 0,$ $\phi''(r)\geq 0$
for all $r\in \left[ 0,R\right] $ and $\phi \left( R\right)=0.$ We refer to \cite{BC} for a proof.

\section{Bottom spectrum estimates}

In this section, we consider upper bounds for the bottom spectrum
of complete stable minimal hypersurfaces.
The bottom spectrum of the Laplacian on a complete manifold $N,$ denoted by $%
\lambda_0(N),$ is an important geometric invariant and characterized as the
optimal Poincar\'e inequality constant or

\begin{equation*}
\lambda _{0}\left( N \right) =\inf_{\phi \in C_{0}^{\infty }\left(N\right) }%
\frac{\int_{N}\left\vert \nabla \phi \right\vert ^{2}}{\int_{N}\phi ^{2}}.
\end{equation*}
According to \cite{LW}, for any $p\in N,$

\begin{equation*}
\lambda _{0}\left( N \right)\leq \frac{1}{4}\,\left(\liminf_{R\to \infty}\,%
\frac{\ln V_p(R)}{R} \right)^2,
\end{equation*}
where $V_p(R)$ denotes the volume of the geodesic ball $B_p(R)$ centered at
point $p$ of radius $R.$ As an
immediate corollary to the area estimate Corollary \ref{i4}, one obtains the
following result which is due to B\'erard, Castillon and Cavalcante \cite{BCC}.

\begin{theorem}[B\'erard, Castillon and Cavalcante]
\label{lambda}Let $\Sigma $ be a complete noncompact stable minimal surface
in a three dimensional manifold $M.$

(a) If the scalar curvature $S$ of $M$ satisfies $S\geq -6$, then 
\begin{equation*}
\lambda _{0}\left( \Sigma \right) \leq 1.
\end{equation*}

(b) If the sectional curvature $K$ of $M$ satisfies $K\geq -1$, then 
\begin{equation*}
\lambda _{0}\left( \Sigma \right) \leq \frac{4}{7}.
\end{equation*}
\end{theorem}

We now provide a different argument for this result.
Recall that an $n-$dimensional manifold $\Sigma $ is called nonparabolic 
if it admits a positive symmetric Green's function. It is well-known that
this is the case if $\lambda_0(\Sigma)>0$ (see \cite{L}).
Let $G(p, x)$ be the minimal positive 
Green's function. Then $G\left( p, x\right)=G\left( x,p\right) >0,$

\begin{equation*}
\Delta _{x}G\left( p,x\right) =-\delta \left( p,x\right)
\end{equation*}%
and 
\begin{equation*}
G\left( p,x\right) =\lim_{i\rightarrow \infty }G_{i}\left( p,x\right),
\end{equation*}%
where $G_{i}\left( p, x\right) $ is the Dirichlet Green's function of a
compact exhaustion $\Omega _{i}$ of $\Sigma.$ 
For fixed point $p,$ we denote $G\left(x\right) =G\left( p,x\right).$ 
It follows from the construction that $\max_{\partial B_p(r)} G$ is a
strictly decreasing function in $r$ and that 

\begin{equation}
\int_{\Sigma \setminus B_{p}\left( 1\right) }\left\vert \nabla G\right\vert
^{2}<\infty.  \label{k}
\end{equation}

Since $G$ is harmonic away from the pole $p,$ the Kato inequality implies

\begin{equation*}
\left\vert G_{ij}\right\vert^{2}\geq \frac{n}{n-1}\left\vert \nabla \left\vert \nabla G\right\vert\right\vert ^{2}.
\end{equation*}
By the Bochner formula we have

\begin{equation}
\Delta \left\vert \nabla G\right\vert \geq \frac{1}{n-1}\,\left\vert \nabla \left\vert \nabla G\right\vert
\right\vert ^{2}\,\left\vert \nabla G\right\vert^{-1}+\mathrm{Ric}^{\Sigma }\left( \nabla G,\nabla G\right) 
\left\vert \nabla G\right\vert^{-1} \label {B1}
\end{equation}
on $\Sigma\setminus \{p\}.$ Similarly, for $v=\ln G,$ using

\begin{equation*}
\Delta v=-\left\vert \nabla v\right\vert^{2}
\end{equation*}
and

\begin{equation*}
\left\vert \nabla v\right\vert\,v_{11}=\left\langle \nabla \left\vert \nabla v\right\vert ,\nabla v\right\rangle,
\end{equation*}
where $\{e_1,\cdots, e_n\}$ is a local orthonormal frame on $\Sigma$ with $e_1=\frac{\nabla}{\left\vert \nabla v\right\vert},$
one easily sees from the Bochner formula together with the standard manipulation that

\begin{eqnarray}
\frac{1}{2}\Delta \left\vert \nabla v\right\vert^2 &\geq& 
\frac{n}{n-1}\,\left\vert \nabla \left\vert \nabla v\right\vert\right\vert ^{2}+\frac{1}{n-1} \left\vert \nabla v\right\vert^4 \label {B2}\\
&&-\frac{n-2}{n-1}\left\langle \nabla \left\vert \nabla v\right\vert^2 ,\nabla v\right\rangle
+\mathrm{Ric}^{\Sigma }\left( \nabla v,\nabla v\right) \notag
\end{eqnarray}
on $\Sigma\setminus \{p\}.$

Denote with 

\begin{eqnarray*}
L\left( a,b\right) &=&\left\{ x\in \Sigma: a<G\left( x\right) <b\right\} \\
l\left( t\right) &=&\left\{ x\in \Sigma: G\left( x\right) =t\right\}.
\end{eqnarray*}
Then $L(\alpha, \infty)\subset B_p(1)$ for $\alpha=\max_{\partial B_p(1)} G.$ 
According to Lemma 5.1 in \cite{LW1},

\begin{eqnarray}
\int_{l\left( t\right) }\left\vert \nabla G\right\vert &=&1  \label{k1} \\
\int_{L\left( a,b\right) }\left\vert \nabla G\right\vert ^{2} f\left(
G\right) &=&\int_{a}^{b} f(t) dt  \notag
\end{eqnarray}%
for any integrable function $f$ provided that $\lambda_0(\Sigma)>0.$ 

We need the following integral gradient estimate. 

\begin{lemma}
\label{G} Let $M$ be a three dimensional manifold with scalar curvature
bounded below. For a complete stable minimal surface $\Sigma $ in $M$ with
$\lambda_0(\Sigma)>0,$
its minimal positive Green's function $G$ satisfies

\begin{equation*}
\int_{\Sigma \setminus B_{p}\left( 1\right) }\frac{\left\vert \nabla
G\right\vert ^{4}}{G^{3}\ln ^{2q}\left( 1+G^{-1}\right) }<\infty 
\end{equation*}
for any $q>\frac{1}{2}.$
\end{lemma}

\begin{proof}
Let $v=\ln G.$ Then by (\ref{B2}),

\begin{equation*}
\frac{1}{2}\Delta \left\vert \nabla v\right\vert ^{2}\ge \left\vert \nabla v\right\vert ^{4}
+2\left\vert \nabla \left\vert \nabla v\right\vert \right\vert^{2}+K_{\Sigma }\left\vert \nabla v\right\vert ^{2}.
\end{equation*}
on $\Sigma\setminus \{p\}.$ According to (\ref{a2}), 

\begin{eqnarray*}
K_{\Sigma }\left\vert \nabla v\right\vert ^{2} &=&\frac{1}{2}S\left\vert
\nabla v\right\vert ^{2}-\left( \mathrm{Ric}\left( \nu ,\nu \right) +\frac{1%
}{2}\left\vert h\right\vert ^{2}\right) \left\vert \nabla v\right\vert ^{2}
\\
&\geq &-C\left\vert \nabla v\right\vert ^{2}-\left( \mathrm{Ric}\left( \nu
,\nu \right) +\left\vert h\right\vert ^{2}\right) \left\vert \nabla
v\right\vert ^{2}.
\end{eqnarray*}%
Thus, for any cut-off function $\phi,$ we have

\begin{eqnarray}
\int_{\Sigma }\left\vert \nabla v\right\vert^{4}\,\phi ^{2} 
&\leq& \frac{1}{2}\int_{\Sigma }\phi ^{2}\,\Delta \left\vert \nabla v\right\vert ^{2}
+C\int_{\Sigma }\phi ^{2}\,\left\vert\nabla v\right\vert ^{2}  \label{k3} \\
&&-2\int_{\Sigma }\phi ^{2}\,\left\vert \nabla \left\vert \nabla v\right\vert \right\vert^{2} \notag \\
&&+\int_{\Sigma }\left( \mathrm{Ric}\left( \nu ,\nu \right) +\left\vert
h\right\vert ^{2}\right) \left\vert \nabla v\right\vert ^{2}\,\phi ^{2}. 
\notag
\end{eqnarray}%
By the stability inequality (\ref{S}), the last term is estimated as

\begin{eqnarray*}
&&\int_{\Sigma }\left( \mathrm{Ric}\left( \nu ,\nu \right) +\left\vert
h\right\vert ^{2}\right) \left\vert \nabla v\right\vert ^{2}\,\phi ^{2} \\
&\leq&\int_{\Sigma }\left\vert \nabla \left( \left\vert \nabla v\right\vert
\phi \right) \right\vert ^{2} \\
&=&\int_{\Sigma }\left\vert \nabla \phi \right\vert ^{2}\left\vert \nabla v\right\vert ^{2}
-\int_{\Sigma }\phi ^{2}\,\left\vert \nabla
v\right\vert\,\Delta \left\vert \nabla v\right\vert.
\end{eqnarray*}
Hence, (\ref{k3}) becomes

\begin{equation*}
\int_{\Sigma }\left\vert \nabla v\right\vert^{4}\,\phi ^{2} \leq
\int_{\Sigma }\left(C\,\phi ^{2}+\left\vert \nabla \phi \right\vert ^{2}\right)\,\left\vert\nabla v\right\vert ^{2}.
\end{equation*}
For $\frac{1}{2}<q<1,$ let $\phi=\psi\,f(G),$ where $\psi$ is a cut-off function
such that $\psi=0$ on $B_p(1)\cup \left(M\setminus B_p(2R)\right),$ $\psi=1$ on $B_p(R)\setminus B_p(2),$
and
 
\begin{equation*}
f\left( G\right) =\frac{G^{\frac{1}{2}}}{\ln ^{q}\left( A\,G^{-1}\right) }
\end{equation*}
with $A=e^{4}\alpha,$ $\alpha=\max_{\partial B_p(1)} G.$
Direct calculations imply

\begin{equation*}
\int_{\Sigma }\left\vert \nabla \phi \right\vert ^{2}\,\left\vert\nabla v\right\vert ^{2}
\leq 4\int_{\Sigma }\left\vert \nabla \psi \right\vert ^{2}\,f^2\,\left\vert\nabla v\right\vert ^{2}
+\frac{3}{4}\,\int_{\Sigma }\left\vert \nabla v\right\vert^{4}\,\phi ^{2}.
\end{equation*}
Therefore, 

\begin{eqnarray*}
\int_{\Sigma }\left\vert \nabla v\right\vert^{4}\,\phi ^{2} &\leq&
C\,\int_{\Sigma }\left(\psi ^{2}+\left\vert \nabla \psi \right\vert ^{2}\right)\,f^2\,\left\vert\nabla v\right\vert ^{2}\\
&\leq& C\,\int_{\Sigma\setminus B_p(1)}\frac{\left\vert \nabla G\right\vert
^{2}}{G\,\left(\ln \frac{A}{G}\right)^{2q}}\\
&\leq& C\,\int_{L(0, \alpha)}\frac{\left\vert \nabla G\right\vert^{2}}{G\,\left(\ln \frac{A}{G}\right)^{2q}}\\
&=&C\,\int_{0}^{\alpha} \frac{1}{t\,\left(\ln \frac{A}{t}\right)^{2q}}\,dt\\
&\leq& C.
\end{eqnarray*}%
Letting $R\to \infty,$ we conclude
 
\begin{equation}
\int_{\Sigma\setminus B_p(2)}\frac{\left\vert \nabla G\right\vert ^{4}}{G^{3}\ln^{2q}\left( AG^{-1}\right) }\leq C.  \label{k7}
\end{equation}%
This proves the result.
\end{proof}

Note that by the Cauchy-Schwarz inequality it follows that 
\begin{eqnarray*}
\left( \int_{\Sigma \backslash B_{p}\left( 1\right) }\frac{\left\vert \nabla
G\right\vert ^{3}}{G^{2}\ln ^{2q}\left( 1+G^{-1}\right) }\right) ^{2} &\leq
&\left( \int_{\Sigma \backslash B_{p}\left( 1\right) }\frac{\left\vert
\nabla G\right\vert ^{4}}{G^{3}\ln ^{2q}\left( 1+G^{-1}\right) }\right) \\
&&\times \left( \int_{\Sigma \backslash B_{p}\left( 1\right) }\frac{%
\left\vert \nabla G\right\vert ^{2}}{G\ln ^{2q}\left( 1+G^{-1}\right) }%
\right) \\
&\leq &C\left( \int_{\Sigma \backslash B_{p}\left( 1\right) }\frac{%
\left\vert \nabla G\right\vert ^{2}}{G\ln ^{2q}\left( 1+G^{-1}\right) }%
\right) .
\end{eqnarray*}%
Again, the last integral is finite for $q>\frac{1}{2}$ by (\ref{k1}) and the co-area formula. 
Therefore,
 
\begin{equation}
\int_{\Sigma \backslash B_{p}\left( 1\right) }\frac{\left\vert \nabla
G\right\vert ^{3}}{G^{2}\ln ^{2q}\left( 1+G^{-1}\right) }<\infty.
\label{k8}
\end{equation}

We are now ready to prove Theorem \ref{lambda}.

\begin{proof} [Proof of Theorem \ref{lambda}]
Without loss of generality we may assume that $\lambda _{0}\left( \Sigma
\right) >0.$ Then $\Sigma $ is nonparabolic. Let $G\left(
p,x\right) $ be the minimal positive Green's function of $\Sigma $ with a
pole at $p\in \Sigma.$ For simplicity, we denote this function by $G\left(
x\right).$ Note that $G$ is harmonic on $\Sigma \backslash \left\{p\right\}.$

Let 
\begin{equation*}
L\left( a,b\right) =\left\{ x\in \Sigma: a<G\left( x\right) <b\right\}.
\end{equation*}%
For $\varepsilon>0$ small enough, define function $\chi $ by

\begin{equation}
\chi \left( x\right) =\left\{ 
\begin{array}{c}
1 \\ 
\frac{\ln G\left( x\right) -\ln \left(\varepsilon^2 \right) }{-\ln
\varepsilon} \\ 
0%
\end{array}%
\right. 
\begin{array}{l}
\text{on }L\left( \varepsilon ,\infty\right) \\ 
\text{on }L\left(\varepsilon^2 ,\varepsilon \right) \\ 
\text{on }L\left(0, \varepsilon\right) \\ 
\end{array}
\label{b1}
\end{equation}%
Since $G\left( x\right) $ may not converge to zero as $x\rightarrow \infty,$
the function $\chi $ may not have compact support. So we consider the
cut-off function

\begin{equation*}
\varphi =\chi \psi,
\end{equation*}
where 
\begin{equation}
\psi \left( x\right) =\left\{ 
\begin{array}{c}
0\\
r(x)-1\\
1 \\ 
R+1-r\left( x\right) \\ 
0%
\end{array}%
\right. 
\begin{array}{l}
\text{on }B_p(1)\\
\text{on }B_p(2)\setminus B_p(1)\\
\text{on }B_{p}\left( R\right)\setminus B_p(2)\\ 
\text{on }B_{p}\left( R+1\right) \setminus B_{p}\left( R\right) \\ 
\text{on }\Sigma \setminus B_{p}\left( R+1\right)%
\end{array}
\label{b2}
\end{equation}

Setting 
\begin{equation*}
\phi =\left\vert \nabla G\right\vert ^{\frac{1}{2}}\varphi 
\end{equation*}
in the Poincar\'e inequality, we have
 
\begin{equation*}
\lambda _{0}\left( \Sigma \right) \int_{\Sigma }\left\vert \nabla
G\right\vert \varphi ^{2}\leq \int_{\Sigma }\left\vert \nabla \left(
\left\vert \nabla G\right\vert ^{\frac{1}{2}}\varphi \right) \right\vert
^{2}.
\end{equation*}%
Expanding the right side, we get

\begin{eqnarray}
\lambda _{0}\left( \Sigma \right) \int_{\Sigma }\left\vert \nabla
G\right\vert \varphi ^{2} &\leq &\frac{1}{4}\int_{\Sigma }\left\vert \nabla
\left\vert \nabla G\right\vert \right\vert ^{2}\left\vert \nabla
G\right\vert ^{-1}\varphi ^{2} \label{b5}\\
&&+\int_{\Sigma }\left\vert \nabla G\right\vert
\left\vert \nabla \varphi \right\vert ^{2}   \notag\\
&&+\frac{1}{2}\int_{\Sigma }\left\langle \nabla \varphi ^{2},\nabla
\left\vert \nabla G\right\vert \right\rangle  \notag \\
&\leq& \left(\frac{1}{4}+\delta\right)\int_{\Sigma }\left\vert \nabla
\left\vert \nabla G\right\vert \right\vert ^{2}\left\vert \nabla
G\right\vert ^{-1}\varphi ^{2} \notag\\
&&+C(\delta)\int_{\Sigma }\left\vert \nabla G\right\vert
\left\vert \nabla \varphi \right\vert ^{2} \notag
\end{eqnarray}
for any $\delta>0.$ We now estimate the first term on the right hand side. Note that by 
(\ref{B1}),

\begin{equation}
\Delta \left\vert \nabla G\right\vert \geq \left\vert \nabla \left\vert
\nabla G\right\vert \right\vert ^{2}\left\vert \nabla G\right\vert
^{-1}+K_{\Sigma }\left\vert \nabla G\right\vert  \label{b6}
\end{equation}%
on $\Sigma\setminus \left\{ p\right\} $ whenever $\left\vert \nabla G\right\vert \neq 0.$

In the case that $S\geq -6$, by (\ref{a2}) it follows that

\begin{equation*}
K_{\Sigma }\left\vert \nabla G\right\vert \geq -3\left\vert \nabla
G\right\vert -\left( \mathrm{Ric}\left( \nu ,\nu \right) +\left\vert
h\right\vert ^{2}\right) \left\vert \nabla G\right\vert.
\end{equation*}%
Hence, (\ref{b6}) becomes

\begin{eqnarray}
\Delta \left\vert \nabla G\right\vert &\geq &\left\vert \nabla \left\vert
\nabla G\right\vert \right\vert ^{2}\left\vert \nabla G\right\vert
^{-1}-3\left\vert \nabla G\right\vert  \label{b7} \\
&&-\left( \mathrm{Ric}\left( \nu ,\nu \right) +\left\vert h\right\vert
^{2}\right) \left\vert \nabla G\right\vert.  \notag
\end{eqnarray}%
Integrating by parts we have

\begin{equation*}
\int_{\Sigma}\varphi ^{2}\Delta \left\vert \nabla G\right\vert
=-\int_{\Sigma}\left\langle \nabla \left\vert \nabla G\right\vert ,\nabla \varphi
^{2}\right\rangle .
\end{equation*}%
Therefore, (\ref{b7}) implies

\begin{eqnarray}
\int_{\Sigma }\left\vert \nabla \left\vert \nabla G\right\vert \right\vert
^{2}\left\vert \nabla G\right\vert ^{-1}\varphi ^{2} &\leq &3\int_{\Sigma
}\left\vert \nabla G\right\vert \varphi ^{2}  \label{b8} \\
&&+\int_{\Sigma }\left( \mathrm{Ric}\left( \nu ,\nu \right) +\left\vert
h\right\vert ^{2}\right) \left\vert \nabla G\right\vert \varphi ^{2}  \notag
\\
&&-\int_{\Sigma }\left\langle \nabla \left\vert \nabla G\right\vert ,\nabla
\varphi ^{2}\right\rangle .  \notag
\end{eqnarray}%
Using the stability inequality (\ref{S}) we have that

\begin{eqnarray*}
\int_{\Sigma }\left( \mathrm{Ric}\left( \nu ,\nu \right) +\left\vert
h\right\vert ^{2}\right) \left\vert \nabla G\right\vert \varphi ^{2} &\leq
&\int_{\Sigma }\left\vert \nabla \left( \left\vert \nabla G\right\vert ^{%
\frac{1}{2}}\varphi \right) \right\vert ^{2} \\
&\leq &\frac{1}{4}\int_{\Sigma }\left\vert \nabla \left\vert \nabla
G\right\vert \right\vert ^{2}\left\vert \nabla G\right\vert ^{-1}\varphi ^{2}
\\
&&+\int_{\Sigma }\left\vert \nabla G\right\vert \left\vert \nabla \varphi
\right\vert ^{2} \\
&&+\frac{1}{2}\int_{\Sigma }\left\langle \nabla \varphi ^{2},\nabla
\left\vert \nabla G\right\vert \right\rangle .
\end{eqnarray*}%
Combining with (\ref{b8}) we obtain that 
\begin{eqnarray*}
\int_{\Sigma }\left\vert \nabla \left\vert \nabla G\right\vert
\right\vert ^{2}\left\vert \nabla G\right\vert ^{-1}\varphi ^{2} &\leq&
4\,\int_{\Sigma }\left\vert \nabla G\right\vert \varphi ^{2}+\frac{4}{3}%
\int_{\Sigma }\left\vert \nabla G\right\vert \left\vert \nabla \varphi
\right\vert ^{2}  \\
&&-\frac{2}{3}\int_{\Sigma }\left\langle \nabla \varphi ^{2},\nabla
\left\vert \nabla G\right\vert \right\rangle \\
&\leq& 4\,\int_{\Sigma }\left\vert \nabla G\right\vert \varphi ^{2}+C(\delta)
\int_{\Sigma }\left\vert \nabla G\right\vert \left\vert \nabla \varphi
\right\vert ^{2}\\
&&+\delta\,\int_{\Sigma }\left\vert \nabla \left\vert \nabla G\right\vert
\right\vert ^{2}\left\vert \nabla G\right\vert ^{-1}\varphi ^{2}. 
\end{eqnarray*}
Therefore,

\begin{equation*}
\int_{\Sigma }\left\vert \nabla \left\vert \nabla G\right\vert
\right\vert ^{2}\left\vert \nabla G\right\vert ^{-1}\varphi ^{2} \leq
\frac{4}{1-\delta}\,\int_{\Sigma }\left\vert \nabla G\right\vert \varphi ^{2}+C(\delta)
\int_{\Sigma }\left\vert \nabla G\right\vert \left\vert \nabla \varphi
\right\vert ^{2}.
\end{equation*}
Plugging into (\ref{b5}) then yields that

\begin{equation}
\left(\lambda _{0}\left( \Sigma \right)-  \frac{1+4\delta}{1-\delta}\right)\int_{\Sigma }\left\vert \nabla
G\right\vert \varphi ^{2} 
\leq C(\delta)\,\int_{\Sigma }\left\vert \nabla
G\right\vert \left\vert \nabla \varphi \right\vert ^{2}.  \label{o2}
\end{equation}
We now estimate the right hand side. Obviously,

\begin{equation}
\int_{\Sigma }\left\vert \nabla G\right\vert \left\vert \nabla \varphi \right\vert ^{2}
\leq 2\int_{\Sigma
}\left\vert \nabla G\right\vert \left\vert \nabla \chi \right\vert ^{2}\psi
^{2}+2\int_{\Sigma }\left\vert \nabla G\right\vert \left\vert \nabla \psi
\right\vert ^{2}\chi ^{2}.  \label{o5}
\end{equation}
The first term of the right hand side is bounded by

\begin{eqnarray*}
\int_{\Sigma}\left\vert \nabla G\right\vert \left\vert \nabla \chi \right\vert ^{2}\psi^{2}
&\leq& \frac{1}{\left(\ln \varepsilon\right)^2}\,
\int_{L(\varepsilon^2, \varepsilon)} \frac{\left\vert \nabla G\right\vert^3}{G^2}\\
&\leq& C\,\frac{1}{\left(\ln \varepsilon\right)^2}\, \ln ^{2q}\left( 1+\varepsilon^{-2}\right)\\
&\leq& C\, \ln ^{2q-2}\left( 1+\varepsilon^{-1}\right),
\end{eqnarray*}
where we have used (\ref{k8}) in the second line.
The second term of (\ref{o5}) is estimated as 

\begin{equation*}
\int_{\Sigma }\left\vert \nabla G\right\vert \left\vert \nabla \psi
\right\vert ^{2}\chi ^{2} \leq \left( \int_{B_p(R+1)\setminus B_p(R) }\left\vert \nabla
G\right\vert ^{2}\right) ^{\frac{1}{2}}
\left( \int_{\left(B_p(R+1)\setminus B_{p}(R)\right)\cap L(\varepsilon^2, \infty)} \chi^2\right)^{\frac{1}{2}}+C.
\end{equation*}
However, the integral estimate in \cite{LW1} says that

\begin{equation*}
\int_{\Sigma \backslash B_{p}\left( R\right) }G^{2}\leq C\,e^{-2\sqrt{\lambda
_{0}\left( \Sigma \right) }R}
\end{equation*}
and
\begin{equation*}
\int_{\Sigma \backslash B_{p}\left( R\right) }\left\vert \nabla G\right\vert^{2}\leq C\,e^{-2\sqrt{\lambda_{0}\left( \Sigma \right) }R}.
\end{equation*}
In particular,

\begin{eqnarray*}
\int_{\left(B_p(R+1)\setminus B_{p}(R)\right)\cap L(\varepsilon^2, \infty)} \chi^2
&\leq &\frac{1}{\varepsilon ^{4}}\int_{\Sigma \backslash B_{p}\left( R\right) }G^{2} \\
&\leq &\frac{C}{\varepsilon ^{4}}e^{-2\sqrt{\lambda_{0}(\Sigma)}\,R}.
\end{eqnarray*}%
Hence, 

\begin{equation*}
\int_{\Sigma }\left\vert \nabla G\right\vert \left\vert \nabla \psi
\right\vert ^{2}\chi ^{2} \leq \frac{C}{\varepsilon ^{2}}e^{-2\sqrt{\lambda_{0}(\Sigma)}\,R}+C.
\end{equation*}
In conclusion, (\ref{o5}) becomes

\begin{equation*}
\int_{\Sigma }\left\vert \nabla G\right\vert \left\vert \nabla \varphi \right\vert ^{2}\leq C\, \ln ^{2q-2}\left( 1+\varepsilon^{-1}\right)
+\frac{C}{\varepsilon ^{2}}e^{-2\sqrt{\lambda_{0}(\Sigma)}\,R}+C.
\end{equation*}
Plugging into (\ref{o2}), we arrive at

\begin{equation*}
\left(\lambda _{0}\left( \Sigma \right)-\frac{1+4\delta}{1-\delta}\right)
\int_{\Sigma }\left\vert \nabla G\right\vert \varphi ^{2} 
\leq C(\delta)\,\left(\ln ^{2q-2}\left( 1+\varepsilon^{-1}\right)
+\varepsilon ^{-2}\,e^{-2\sqrt{\lambda_{0}(\Sigma)}\,R}+C\right)
\end{equation*}
By first letting $R\to \infty$ and then $\varepsilon\to 0,$ as $q<1,$
we conclude that

\begin{equation}
\left(\lambda _{0}\left( \Sigma \right)-\frac{1+4\delta}{1-\delta}\right)
\int_{\Sigma\setminus B_p(2) }\left\vert \nabla G\right\vert
\leq C(\delta). \label{o6}
\end{equation}
However, using

\begin{equation*}
\int_{\partial B_p(r)} \frac{\partial G}{\partial r}=-1,
\end{equation*}
we have

\begin{equation*}
\int_{\Sigma\setminus B_p(2) }\left\vert \nabla G\right\vert
\geq -\int_2^{\infty}dr \int_{\partial B_p(r)} \frac{\partial G}{\partial r}=\infty.
\end{equation*}
By (\ref{o6}), it follows that

\begin{equation*}
\lambda _{0}\left( \Sigma \right)\leq \frac{1+4\delta}{1-\delta}
\end{equation*}
for any $\delta>0.$
Therefore, $\lambda _{0}\left( \Sigma \right)\leq 1.$

Now we consider the case that $K\geq -1.$ By (\ref{a1}) we have

\begin{equation*}
K_{\Sigma }\left\vert \nabla G\right\vert \geq -\left\vert \nabla
G\right\vert -\frac{1}{2}\left\vert h\right\vert ^{2}\left\vert \nabla
G\right\vert.
\end{equation*}
Therefore, (\ref{b6}) becomes

\begin{equation}
\Delta \left\vert \nabla G\right\vert \geq \left\vert \nabla \left\vert
\nabla G\right\vert \right\vert ^{2}\left\vert \nabla G\right\vert
^{-1}-\left\vert \nabla G\right\vert -\frac{1}{2}\left\vert h\right\vert
^{2}\left\vert \nabla G\right\vert .  \label{b10}
\end{equation}%
However, by the stability inequality (\ref{S}),

\begin{eqnarray*}
\int_{\Sigma }\left( \left\vert h\right\vert ^{2}+\mathrm{Ric}\left( \nu
,\nu \right) \right) \left\vert \nabla G\right\vert \varphi ^{2} &\leq &%
\frac{1}{4}\int_{\Sigma }\left\vert \nabla \left\vert \nabla G\right\vert
\right\vert ^{2}\left\vert \nabla G\right\vert ^{-1}\varphi ^{2} \\
&&+\frac{1}{2}\int_{\Sigma }\left\langle \nabla \left\vert \nabla
G\right\vert ,\nabla \varphi ^{2}\right\rangle \\
&&+\int_{\Sigma }\left\vert \nabla G\right\vert \left\vert \nabla \varphi
\right\vert ^{2}.
\end{eqnarray*}%
Hence, using that $\mathrm{Ric}\left( \nu ,\nu \right) \geq -2$, we conclude

\begin{eqnarray}
\frac{1}{2}\int_{\Sigma }\left\vert h\right\vert ^{2}\left\vert \nabla
G\right\vert \varphi ^{2} &\leq &\frac{1}{8}\int_{\Sigma }\left\vert \nabla
\left\vert \nabla G\right\vert \right\vert ^{2}\left\vert \nabla
G\right\vert ^{-1}\varphi ^{2}+\int_{\Sigma }\left\vert \nabla G\right\vert
\varphi ^{2}  \label{b11} \\
&&+\frac{1}{4}\int_{\Sigma }\left\langle \nabla \left\vert \nabla
G\right\vert ,\nabla \varphi ^{2}\right\rangle +\frac{1}{2}\int_{\Sigma
}\left\vert \nabla G\right\vert \left\vert \nabla \varphi \right\vert ^{2}. 
\notag
\end{eqnarray}%
By (\ref{b10}) and (\ref{b11}) it follows that for any small $\delta>0,$

\begin{eqnarray*}
\int_{\Sigma }\left\vert \nabla \left\vert \nabla G\right\vert \right\vert
^{2}\left\vert \nabla G\right\vert ^{-1}\varphi ^{2} &\leq &\frac{1}{2}%
\int_{\Sigma }\left\vert h\right\vert ^{2}\left\vert \nabla G\right\vert
\varphi ^{2}+\int_{\Sigma }\left\vert \nabla G\right\vert \varphi
^{2}-\int_{\Sigma }\left\langle \nabla \left\vert \nabla G\right\vert
,\nabla \varphi ^{2}\right\rangle \\
&\leq &\left(\frac{1}{8}+\delta\right)\int_{\Sigma }\left\vert \nabla \left\vert \nabla
G\right\vert \right\vert ^{2}\left\vert \nabla G\right\vert ^{-1}\varphi ^{2}
+2\,\int_{\Sigma }\left\vert \nabla G\right\vert \varphi^{2}\\
&&+C(\delta)\,\int_{\Sigma}\left\vert \nabla G\right\vert \left\vert \nabla \varphi \right\vert ^{2}.
\end{eqnarray*}%
In conclusion,

\begin{equation*}
\int_{\Sigma }\left\vert \nabla \left\vert \nabla G\right\vert \right\vert
^{2}\left\vert \nabla G\right\vert ^{-1}\varphi ^{2} \leq \frac{16}{7-8\delta}
\int_{\Sigma }\left\vert \nabla G\right\vert \varphi ^{2}+C(\delta)\,\int_{\Sigma }\left\vert \nabla G\right\vert \left\vert \nabla \varphi
\right\vert ^{2}.
\end{equation*}
Plugging into (\ref{b5}), we have

\begin{equation*}
\left(\lambda _{0}\left( \Sigma \right)- \frac{4+16\delta}{7-8\delta}\right)
\int_{\Sigma }\left\vert \nabla
G\right\vert \varphi ^{2} \leq 
C(\delta)\,\int_{\Sigma }\left\vert \nabla G\right\vert \left\vert \nabla \varphi
\right\vert ^{2}.
\end{equation*}
Similarly, we can conclude that

\begin{equation*}
\lambda _{0}\left( \Sigma \right)\leq \frac{4+16\delta}{7-8\delta}.
\end{equation*}
Now letting $\delta\to 0,$ we have $\lambda _{0}\left( \Sigma \right)\leq \frac{4}{7}.$
\end{proof}

The preceding argument can be generalized to stable minimal hypersurfaces of
dimension up to five. Let $\left( M,g\right) $ be an $(n+1)$-dimensional
complete Riemannian manifold with its sectional curvature bounded below by 

\begin{equation*}
K\geq -1.
\end{equation*}
Let $\Sigma \subset M$ be a stable minimal hypersurface in $M.$ 
Then the stability inequality (\ref{S}) implies that 

\begin{equation}
\int_{\Sigma }\left\vert h\right\vert ^{2}\phi ^{2}\leq \int_{\Sigma
}\left\vert \nabla \phi \right\vert ^{2}+n\,\int_{\Sigma }\phi ^{2}.
\label{z1}
\end{equation}%
For a local orthonormal frame $\left\{ e_{1},\cdots, e_{n}\right\}$ 
of $\Sigma,$ by the Gauss curvature equations,

\begin{eqnarray}
R_{aa}^{\Sigma } &=&\sum_{c}R_{acac}-\sum_{c}\left\vert h_{ac}\right\vert
^{2}  \label{z2} \\
&\geq &-(n-1)-\frac{n-1}{n}\left\vert h\right\vert ^{2}  \notag
\end{eqnarray}
for indices $1\leq a, \,c \leq n,$ where in the last line we have used that $\Sigma $ is minimal. The argument in
Lemma \ref{G} can be carried over to prove the following.

\begin{lemma}
\label{G'} Let $M$ be an $(n+1)$-dimensional complete manifold with sectional curvature
bounded below and $n\leq 5.$ For a complete stable minimal hypersurface $\Sigma$ in $M$ with
$\lambda_0(\Sigma)>0,$
its minimal positive Green's function $G\left( x\right) =G\left(p,x\right) $ satisfies

\begin{equation*}
\int_{\Sigma \setminus B_{p}\left( 1\right) }
\frac{\left\vert \nabla G\right\vert ^{4}}{G^{3}\ln ^{2q}\left( 1+G^{-1}\right) }<\infty
\end{equation*}%
for any $q>\frac{1}{2}.$
\end{lemma}

\begin{proof}
Let $v=\ln G$. Then, according to (\ref{B2}),

\begin{eqnarray*}
\frac{1}{2}\Delta \left\vert \nabla v\right\vert^2 &\geq& \frac{1}{n-1} \left\vert \nabla v\right\vert^4+
\frac{n}{n-1}\,\left\vert \nabla \left\vert \nabla v\right\vert\right\vert ^{2}\\
&&-\frac{n-2}{n-1}\left\langle \nabla \left\vert \nabla v\right\vert^2 ,\nabla v\right\rangle
+\mathrm{Ric}^{\Sigma }\left( \nabla v,\nabla v\right) 
\end{eqnarray*}
on $\Sigma \setminus \{p\}.$ Note that by (\ref{z2}), 

\begin{equation*}
\mathrm{Ric}^{\Sigma}(\nabla v, \nabla v)\geq -C\left\vert \nabla v\right\vert ^{2}-\frac{n-1}{n}\left( \mathrm{Ric}\left( \nu
,\nu \right) +\left\vert h\right\vert ^{2}\right) \left\vert \nabla v\right\vert ^{2}
\end{equation*}
as the sectional curvature of $M$ is bounded from below.
Thus, for any cut-off function $\phi,$ noting that

\begin{equation*}
\int_{\Sigma }\left( \mathrm{Ric}\left( \nu ,\nu \right) +\left\vert
h\right\vert ^{2}\right) \left\vert \nabla v\right\vert ^{2}\,\phi ^{2} 
\leq \int_{\Sigma }\left\vert \nabla \left( \left\vert \nabla v\right\vert
\phi \right) \right\vert ^{2}
\end{equation*}
by the stability inequality (\ref{S}), we have

\begin{eqnarray}
&&\frac{1}{n-1}\int_{\Sigma }\left\vert \nabla v\right\vert^{4}\,\phi ^{2} \label{q3} \\
&\leq& \frac{1}{2}\int_{\Sigma }\phi ^{2}\,\Delta \left\vert \nabla v\right\vert ^{2}
+C\int_{\Sigma }\phi ^{2}\,\left\vert\nabla v\right\vert ^{2}  
-\frac{n}{n-1}\,\int_{\Sigma }\phi ^{2}\,\left\vert \nabla \left\vert \nabla v\right\vert \right\vert^{2} \notag\\
&&+\frac{n-2}{n-1}\int_{\Sigma }\phi ^{2}\, \left\langle \nabla \left\vert \nabla v\right\vert^2,\nabla v\right\rangle 
+\frac{n-1}{n}\int_{\Sigma }\left\vert \nabla \left( \left\vert \nabla v\right\vert
\phi \right) \right\vert ^{2} \notag \\
&=&-\frac{1}{n}\int_{\Sigma }\phi\,\left\langle \nabla \phi,\nabla \left\vert \nabla v\right\vert ^{2}\right\rangle
+C\int_{\Sigma }\phi ^{2}\,\left\vert\nabla v\right\vert ^{2} 
-\frac{2n-1}{n(n-1)}\,\int_{\Sigma }\phi ^{2}\,\left\vert \nabla \left\vert \nabla v\right\vert \right\vert^{2} \notag \\
&&+\frac{n-2}{n-1}\int_{\Sigma }\phi ^{2}\, \left\langle \nabla \left\vert \nabla v\right\vert^2,\nabla v\right\rangle 
+\frac{n-1}{n}\int_{\Sigma }\left\vert \nabla v\right\vert^2\,\left\vert \nabla \phi \right\vert ^{2}. \notag 
\end{eqnarray}
Let $\phi=G^{\frac{1}{2}}\,\eta,$ where $\eta$ is a cut-off function on $\Sigma$ with $\eta=0$ on $B_p(1).$
Note that for $\delta>0,$ 

\begin{equation*}
\int_{\Sigma }\left\vert \nabla v\right\vert^2\,\left\vert \nabla \phi \right\vert ^{2}
\leq C(\delta)\,\int_{\Sigma }\left\vert \nabla v\right\vert^2\,\left\vert \nabla \eta \right\vert ^{2}\,G
+\frac{1+\delta}{4}\int_{\Sigma }\phi^2\,\left\vert \nabla v\right\vert^4.
\end{equation*}
Therefore, we conclude from (\ref{q3}) that

\begin{eqnarray}
&&\left(\frac{1}{n-1}-\frac{1+\delta}{4}\,\frac{n-1}{n}\right)
\int_{\Sigma }\left\vert \nabla v\right\vert^{4}\,\phi ^{2} \label{q4}\\
&\leq& C(\delta)\,\int_{\Sigma }\left(\eta ^{2}+\left\vert \nabla \eta\right\vert ^{2}\right)G\,\left\vert\nabla v\right\vert ^{2} 
-\frac{1}{2n}\int_{\Sigma }\left\langle \nabla (G\eta^2),\nabla \left\vert \nabla v\right\vert ^{2}\right\rangle \notag\\
&&-\frac{2n-1}{n(n-1)}\,\int_{\Sigma }\phi ^{2}\,\left\vert \nabla \left\vert \nabla v\right\vert \right\vert^{2} 
+\frac{n-2}{n-1}\int_{\Sigma }\eta^{2}\, \left\langle \nabla \left\vert \nabla v\right\vert^2,\nabla G\right\rangle \notag\\
&=&C(\delta)\,\int_{\Sigma }\left(\eta ^{2}+\left\vert \nabla \eta\right\vert ^{2}\right)G\,\left\vert\nabla v\right\vert ^{2} 
-\frac{1}{2n}\int_{\Sigma }\left\langle \nabla \eta^2,\nabla \left\vert \nabla v\right\vert ^{2}\right\rangle\,G\notag\\
&&-\frac{2n-1}{n(n-1)}\,\int_{\Sigma }\phi ^{2}\,\left\vert \nabla \left\vert \nabla v\right\vert \right\vert^{2} 
+\left(\frac{n-2}{n-1}-\frac{1}{2n}\right)\int_{\Sigma }\eta^{2}\, \left\langle \nabla \left\vert \nabla v\right\vert^2,\nabla G\right\rangle. \notag
\end{eqnarray}
However, as $G$ is harmonic on $\Sigma\setminus \{p\},$

\begin{equation*}
\int_{\Sigma }\eta^{2}\, \left\langle \nabla \left\vert \nabla v\right\vert^2,\nabla G\right\rangle
=-\int_{\Sigma }\left\langle \nabla \eta^2,\nabla G \right\rangle\,\left\vert \nabla v\right\vert ^{2}.
\end{equation*}
Plugging into (\ref{q4}), we get

\begin{eqnarray}
&&\left(\frac{1}{n-1}-\frac{1+\delta}{4}\,\frac{n-1}{n}\right)
\int_{\Sigma }\left\vert \nabla v\right\vert^{4}\,\phi ^{2} \label{q5}\\
&\leq&C(\delta)\,\int_{\Sigma }\left(\eta ^{2}+\left\vert \nabla \eta\right\vert ^{2}\right)G\,\left\vert\nabla v\right\vert ^{2} 
-\frac{2}{n}\int_{\Sigma }\left\langle \nabla \eta,\nabla \left\vert \nabla v\right\vert\right\rangle\,G\eta
\left\vert \nabla v\right\vert \notag\\
&&-\frac{2n-1}{n(n-1)}\,\int_{\Sigma }\eta^{2} G\,\left\vert \nabla \left\vert \nabla v\right\vert \right\vert^{2} 
-\left(\frac{n-2}{n-1}-\frac{1}{2n}\right)\int_{\Sigma }\left\langle \nabla \eta^{2},\nabla G\right\rangle \left\vert \nabla v\right\vert^2. \notag\\
&\leq&C(\delta)\,\int_{\Sigma }\left(\eta ^{2}+\left\vert \nabla \eta\right\vert ^{2}\right)G\,\left\vert\nabla v\right\vert ^{2}  
+\delta\,\int_{\Sigma }\left\vert \nabla v\right\vert^{4}\,\phi ^{2}. \notag
\end{eqnarray}
Therefore, 

\begin{equation*}
\left(\frac{1}{n-1}-\frac{1+\delta}{4}\,\frac{n-1}{n}-\delta\right)
\int_{\Sigma }\left\vert \nabla v\right\vert^{4}\,\phi ^{2} 
\leq
C(\delta)\,\int_{\Sigma }\left(\eta ^{2}+\left\vert \nabla \eta\right\vert ^{2}\right)G\,\left\vert\nabla v\right\vert ^{2}. 
\end{equation*}
Since $n\leq 5,$ one may choose $\delta=\delta (n)>0$ such that
$\frac{1}{n-1}-\frac{1+\delta}{4}\,\frac{n-1}{n}-\delta>0.$ In conclusion, there exists an absolute constant $\Gamma>0$
such that

\begin{equation}
\int_{\Sigma }\left\vert \nabla v\right\vert^{4}\,\phi ^{2} 
\leq \Gamma\,\int_{\Sigma }\left(\eta ^{2}+\left\vert \nabla \eta\right\vert ^{2}\right)G\,\left\vert\nabla v\right\vert ^{2} \label{q6}
\end{equation}
for any cut-off function $\eta$ satisfying $\eta=0$ on $B_p(1).$

For $\frac{1}{2}<q<1,$ let $\eta=\psi\,w(G),$ where $\psi$ is a cut-off function
such that $\psi=0$ on $B_p(1)\cup \left(M\setminus B_p(2R)\right),$ $\psi=1$ on $B_p(R)\setminus B_p(2),$
and
 
\begin{equation*}
w\left( G\right) =\frac{1}{\ln ^{q}\left( A\,G^{-1}\right) }
\end{equation*}
with $A=e^{2\sqrt{\Gamma}}\alpha,$ $\alpha=\max_{\partial B_p(1)} G.$
Direct calculations imply

\begin{equation*}
\int_{\Sigma }\eta ^{2}\,G\,\left\vert\nabla v\right\vert ^{2}\leq \int_{L(0, \alpha)} 
\frac{\left\vert\nabla G\right\vert ^{2}}{G\ln^{2q}\left( A\,G^{-1}\right)}\leq C
\end{equation*}
and

\begin{eqnarray*}
\int_{\Sigma }\left\vert \nabla \eta\right\vert ^{2}\,G\,\left\vert\nabla v\right\vert ^{2}
&\leq& 2\int_{\Sigma }\left\vert \nabla \psi\right\vert ^{2}\,\frac{\left\vert\nabla G\right\vert ^{2}}{G\ln^{2q}\left( A\,G^{-1}\right)}\\
&&+2\int_{\Sigma }\psi^2\left\vert \nabla w\right\vert ^{2}\,G\,\left\vert\nabla v\right\vert ^{2}\\
&\leq& C+\frac{1}{2\Gamma}\,\int_{\Sigma }\phi^2\,\left\vert \nabla v\right\vert ^{4}.
\end{eqnarray*}
Together with (\ref{q6}), we arrive at

\begin{equation*}
\int_{\Sigma }\left\vert \nabla v\right\vert^{4}\,\phi ^{2} 
\leq C+\frac{1}{2}\,\int_{\Sigma }\phi^2\,\left\vert \nabla v\right\vert ^{4}.
\end{equation*}
In other words,

\begin{equation*}
\int_{\Sigma }\left\vert \nabla v\right\vert^{4}\,\phi ^{2} \leq C.
\end{equation*}
Finally, letting $R\to \infty,$ one concludes that

\begin{equation*}
\int_{\Sigma\setminus B_p(2)} \frac{\left\vert\nabla G\right\vert ^{4}}{G^3\ln^{2q}\left( A\,G^{-1}\right)}<\infty.
\end{equation*}
This proves the desired result.
\end{proof}

We are now ready to prove the following spectral estimate.

\begin{theorem}
Let $\Sigma $ be a complete stable minimal hypersurface in $(n+1)$-dimensional
manifold $M$with $n\leq 5.$  If the sectional curvature of $M$ satisfies $K\geq -\kappa$ for some 
nonnegative constant $\kappa,$ then 

\begin{equation*}
\lambda _{0}\left( \Sigma \right) \leq \frac{2n(n-1)^2}{6n-n^2-1}\,\kappa.
\end{equation*}
\end{theorem}

\begin{proof}
Without loss of generality we
assume that $\kappa=1$ and $\lambda _{0}\left( \Sigma \right) >0.$ In particular,
$\Sigma $ is nonparabolic. For fixed $\varepsilon>0$ small enough,
define $\chi $ and $\psi$ by (\ref{b1}) and (\ref{b2}), respectively,
and let

\begin{equation*}
\varphi =\chi \psi.
\end{equation*}%
Setting 
\begin{equation*}
\phi =\left\vert \nabla G\right\vert ^{\frac{1}{2}}\varphi 
\end{equation*}
in the Poincar\'e inequality
and expanding the right side, we get

\begin{equation}
\lambda _{0}\left( \Sigma \right) \int_{\Sigma }\left\vert \nabla
G\right\vert \varphi ^{2} 
\leq \left(\frac{1}{4}+\delta\right)\int_{\Sigma }\left\vert \nabla
\left\vert \nabla G\right\vert \right\vert ^{2}\left\vert \nabla
G\right\vert ^{-1}\varphi ^{2} 
+C(\delta)\int_{\Sigma }\left\vert \nabla G\right\vert
\left\vert \nabla \varphi \right\vert ^{2} \label{q10}
\end{equation}
for any $\delta>0.$ We now estimate the first term on the right hand side. Note that by
(\ref{B1}) and (\ref{z2}),

\begin{eqnarray*}
\Delta \left\vert \nabla G\right\vert &\geq &\frac{1}{n-1}\left\vert \nabla \left\vert \nabla G\right\vert
\right\vert ^{2}\left\vert \nabla G\right\vert^{-1}+\mathrm{Ric}%
^{\Sigma }\left( \nabla G,\nabla G\right) \left\vert \nabla G\right\vert^{-1} \\
&\geq &\frac{1}{n-1}\left\vert \nabla \left\vert \nabla G\right\vert
\right\vert ^{2}\left\vert \nabla G\right\vert^{-1}-
\left( (n-1)+\frac{n-1}{n}\left\vert h\right\vert ^{2}\right) \left\vert \nabla G\right\vert.
\end{eqnarray*}
Hence, 

\begin{eqnarray}
\frac{1}{n-1}\int_{\Sigma }\left\vert \nabla \left\vert \nabla G\right\vert
\right\vert ^{2}\left\vert \nabla G\right\vert ^{-1}\varphi ^{2} &\leq
&(n-1)\int_{\Sigma }\left\vert \nabla G\right\vert \varphi ^{2} \label{b8'} \\
&&+\frac{n-1}{n}\int_{\Sigma }\left\vert h\right\vert ^{2}\left\vert \nabla G\right\vert \varphi ^{2}  \notag\\
&&-\int_{\Sigma }\left\langle \nabla \left\vert \nabla G\right\vert ,\nabla
\varphi ^{2}\right\rangle .  \notag
\end{eqnarray}%
Using the stability inequality (\ref{z1}) we have that

\begin{eqnarray*}
\int_{\Sigma }\left\vert h\right\vert ^{2}\left\vert \nabla G\right\vert
\varphi ^{2} &\leq &\int_{\Sigma }\left\vert \nabla \left( \left\vert \nabla
G\right\vert ^{\frac{1}{2}}\varphi \right) \right\vert ^{2}+n\int_{\Sigma
}\left\vert \nabla G\right\vert \varphi ^{2} \\
&\leq &\frac{1}{4}\int_{\Sigma }\left\vert \nabla \left\vert \nabla
G\right\vert \right\vert ^{2}\left\vert \nabla G\right\vert ^{-1}\varphi
^{2}+\int_{\Sigma }\left\vert \nabla G\right\vert \left\vert \nabla \varphi
\right\vert ^{2} \\
&&+\frac{1}{2}\int_{\Sigma }\left\langle \nabla \varphi ^{2},\nabla
\left\vert \nabla G\right\vert \right\rangle +n\int_{\Sigma }\left\vert
\nabla G\right\vert \varphi ^{2}.
\end{eqnarray*}%
Combining with (\ref{b8'}) we obtain that 

\begin{eqnarray*}
\left(\frac{1}{n-1}-\frac{n-1}{4n}\right)\int_{\Sigma }\left\vert \nabla \left\vert \nabla G\right\vert
\right\vert ^{2}\left\vert \nabla G\right\vert ^{-1}\varphi ^{2} &\leq
&2(n-1)\int_{\Sigma }\left\vert \nabla G\right\vert \varphi ^{2}+C\int_{\Sigma
}\left\vert \nabla G\right\vert \left\vert \nabla \varphi \right\vert ^{2} \\
&&+C\int_{\Sigma }\left\vert \left\langle \nabla \varphi ^{2},\nabla
\left\vert \nabla G\right\vert \right\rangle \right\vert\\
&\leq& 2(n-1)\int_{\Sigma }\left\vert \nabla G\right\vert \varphi ^{2}+C(\delta)\,\int_{\Sigma
}\left\vert \nabla G\right\vert \left\vert \nabla \varphi \right\vert ^{2} \\
&&+\delta\,\int_{\Sigma }\left\vert \nabla \left\vert \nabla G\right\vert
\right\vert ^{2}\left\vert \nabla G\right\vert ^{-1}\varphi ^{2}.
\end{eqnarray*}
Therefore,

\begin{eqnarray*}
\int_{\Sigma }\left\vert \nabla \left\vert \nabla G\right\vert
\right\vert ^{2}\left\vert \nabla G\right\vert ^{-1}\varphi ^{2}&\leq&
\frac{8n(n-1)^2}{6n-n^2-1-4n(n-1)\delta}\int_{\Sigma }\left\vert \nabla G\right\vert \varphi ^{2}\\
&&+C(\delta)\,\int_{\Sigma
}\left\vert \nabla G\right\vert \left\vert \nabla \varphi \right\vert ^{2}.
\end{eqnarray*}
Plugging into (\ref{q10}) then yields that

\begin{equation*}
\left( \lambda _{0}\left( \Sigma \right) -\left(\frac{1}{4}+\delta\right)\frac{8n(n-1)^2}{6n-n^2-1-4n(n-1)\delta}
\right) \int_{\Sigma }\left\vert \nabla G\right\vert \varphi ^{2}\leq C(\delta)\,\int_{\Sigma
}\left\vert \nabla G\right\vert \left\vert \nabla \varphi \right\vert ^{2}.
\end{equation*}
Using Lemma \ref{G'}, one concludes as before that

\begin{equation*}
\lambda_{0}\left( \Sigma \right) \leq \frac{2n(n-1)^2}{6n-n^2-1}.
\end{equation*}
This proves the result.
\end{proof}

\end{document}